\documentclass[11pt]{article}
\title {A DC Programming Approach for Solving Multicast Network Design Problems via the Nesterov Smoothing Technique}
\usepackage{graphicx}
\usepackage{amsmath,amsfonts,amssymb}
\usepackage{url}
\usepackage{amsmath}
\usepackage{epsfig}
\usepackage{epstopdf}
\usepackage{color}
\usepackage{float}
\usepackage{subcaption}
\usepackage{multirow}
\usepackage{algorithm,algpseudocode}
\def\tto{\;{\lower 1pt \hbox{$\rightarrow$}}\kern -10pt
\hbox{\raise 2pt \hbox{$\rightarrow$}}\;}

\def\ra{\rangle}
\def\la{\langle}

\def\epsilon{\varepsilon}
\def\B{\Bbb B}
\def\h{\hfill\Box}
\def\R{\Bbb R}

\def\N{\Bbb N}
\def\ox{\bar{x}}

\def\h{\hfill\square}

\def\ph{\varphi}

\newcounter{lk}

\begin{document}
\begin{center}
{\bf A DC Programming Approach for Solving Multicast Network Design Problems via the Nesterov Smoothing Technique  }\\[2ex]
 W. GEREMEW,\footnote{School of General Studies, Stockton University, Galloway, NJ 08205, USA (wgeremew24@gmail.com). Geremew's research was partly supported by the AFRL Mathematical Modeling and Optimization Institute.} N. M. NAM,\footnote{Fariborz Maseeh Department of Mathematics and Statistics, Portland State University, Portland, OR 97207, USA (mau.nam.nguyen@pdx.edu). Research of this author was partly supported by the National Science Foundation under grant \#1716057.} A. SEMENOV,\footnote{University of Jyv\"askyl\"a, P.O.Box 35 FI-40014 University of Jyvaskyla, Finland (alexander.v.semenov@jyu.fi)} V. BOGINSKI,\footnote{University of Central Florida, 12800 Pegasus Dr., Orlando, FL 32816 USA (vladimir.boginski@ucf.edu)}, and   E. PASILIAO\footnote{Air Force Research Laboratory Eglin AFB, FL, USA. (elpasiliao@gmail.com)}  \\[2ex]

\end{center}
\small{\bf Abstract.} This paper continues our effort initiated in \cite{wondi} to study  \emph{Multicast Communication Networks}, modeled as bilevel hierarchical clustering problems, by using mathematical optimization techniques. Given a finite number of nodes, we consider two different models of  multicast networks by identifying a certain number of nodes as cluster centers, and at the same time, locating a particular node that serves as a total center so as to minimize the total transportation cost through the network. The fact that the cluster centers and the total center have to be among the given nodes makes this problem a discrete optimization problem. Our approach is to reformulate the discrete problem as a continuous one and  to apply Nesterov smoothing approximation technique on the Minkowski gauges that are used as distance measures. This approach enables us to propose two implementable DCA-based algorithms for solving the problems. Numerical results and practical applications are provided to illustrate our approach.\\[1ex]

\newtheorem{Theorem}{Theorem}[section]
\newtheorem{Proposition}[Theorem]{Proposition}
\newtheorem{Remark}[Theorem]{Remark}
\newtheorem{Lemma}[Theorem]{Lemma}
\newtheorem{Corollary}[Theorem]{Corollary}
\newtheorem{Definition}[Theorem]{Definition}
\newtheorem{Example}[Theorem]{Example}
\renewcommand{\theequation}{\thesection.\arabic{equation}}
\normalsize
\vspace*{-0.2in}

\section{Introduction}
\setcounter{equation}{0}

The complexity of modern networks such as communication networks, broadcasting networks, and distribution networks requires  multilevel connectivity. For instance, many department stores usually get their merchandise delivered to them by a delivery company. For efficiency purposes, the delivery company usually wants to identify a certain number of locations  to serve as distribution centers for the delivery of supplies to the stores. At the same time, the company wants to identify a location  as a main distribution center, also  known as the total center, from which the other distribution centers receive their supplies.  This is a typical description of a  bilevel multicast communication network, which can also be seen as a multifacility location problem or as a bilevel hierarchical clustering problem. Borrowing some language from network optimization literature, these problems can be described mathematically as follows:  Given $m$ nodes $a^1, a^2, \ldots, a^m$ in $\R^n$, the objective is to choose $k$ cluster centroids $a^{(1)}$, $a^{(2)}$, \ldots , $a^{(k)}$  and a total center $a^{(k+1)}$ from the given  nodes in such a way that the total \emph{transportation cost} of the tree formed by connecting the cluster centers to the total center, and  the remaining nodes to  the nearest cluster centers is minimized.  The fact that the centers and the total center have to be among the existing nodes makes the problem a discrete optimization problem, which can be shown to be NP--hard.

Many existing algorithms for solving bilevel hierarchical clustering problems are heuristics in nature, and do not optimize any well-defined objective function. The mathematical optimization approach for solving  hierarchical clustering problems was initiated in the pioneering work from \cite{bjor03}. The authors introduced three models of hierarchical clustering based on the Euclidean norm and employed the derivative-free  method developed in \cite{B99} to solve the problem in two dimensions. Replacing the Euclidean norm by the squared Euclidean norm, the authors in \cite{am07} used the DCA, a well-known algorithm for minimizing differences of convex functions introduced by Pham Dinh Tao (see \cite{at97,ta98}), to solve the problem in high dimensions. In fact, the DCA  provides an effective tool for solving the classical clustering problem and its variants; see \cite{abt07,amt14,am07,bjor03,btu16,bvx13,nrg17,wondi} and the references therein. In  our recent work \cite{wondi}, we proposed a new method based on the Nesterov smoothing technique and the DCA to cope with the original models of hierarchical clustering introduced in \cite{bjor03}. The idea of using  the Nesterov smoothing technique overcomes the drawback of the DCA stated in \cite{am07} as ``the DCA is not appropriate for these models''. Our current paper continues the effort initiated in \cite{am07,bjor03} in which mathematical optimization techniques for solving optimization problems beyond convexity are used in multifacility location and clustering. In particular, this paper is the second part of our  paper \cite{wondi} as we propose  other two bivelel  hierarchical clustering models. Another novel component of the present paper compared to \cite{wondi} is the possibility of considering problems with generalized distance generated by Minkowski gauges as well as the possibility to handle problems with constraints.

In this paper, we propose two implementable algorithms based on a DC programming approach combined with the Nesterov smoothing technique to solve the resulting constrained minimization problems for both models. It is important to note that the DCA can only guarantees the convergence to a critical point, so to achieve better results we often run the algorithms
multiple times with different starting points via suitable initialization techniques, such as running the k-means or a genetic algorithm to generate starting centers for the two proposed algorithms.

The paper is organized as follows. In Section 2, we present the continuous optimization formulations of the two models using Minkowski gauges as distance measures. In section 3 we discuss some basic definitions and tools of optimization that are used throughout the paper. In Sections 4 and 5, we develop the two algorithms for the two proposed multicast communication networks. In Section 6 we present our numerical experiments and results performed on artificial datasets as well as real datasets.

\section{Problems Formulation}
\setcounter{equation}{0}

 In this section, we discuss two models of bilevel hierarchical clustering  and provide the tools of optimization  used throughout the paper. In order to reformulate the discrete optimization problem under consideration as a continuous optimization problem, we introduce $k$ artificial centers which are not necessarily  the existing nodes in designing the optimal multicast networks. Denote the $k$ artificial cluster centers by $x^1, x^2, \ldots , x^k$ and the distance measurement between the artificial center $x^\ell ,\;\ell =1, \ldots, k,$ and the real node $a^i, \;i = 1, \ldots, m,$ by a generalized distance $\sigma_F(x^\ell - a^i)$, where $\sigma_F$ is the \emph{support function} associated with a nonempty closed bounded convex set $F$ containing the origin in its interior, i.e.,
 \begin{equation*}
 \sigma_F(x):=\sup\{\la x, y\ra\; |\; y\in F\}.
 \end{equation*}
 Note that if $F$ is the closed unit Euclidean ball in $\R^n$, then $\sigma_F(x)$ defines the Euclidean norm of $x\in \R^n$. In the case where $F$ is the closed unit box of $\R^n$, i.e., $F:=\{u=(u_1, \ldots, u_n)\in \R^n\; |\; -1\leq u_i\leq 1\; \mbox{\rm for }i=1, \ldots, n\}$, then $\sigma_F(x)$ defines the $\ell^1-$norm $\|x\|_1$ of $x\in \R^n$.

 In the first model, the $m$ nodes are clustered around the $k$ artificial centers by trying to minimize the minimum sum of the distances from each node to the $k$ cluster centers. A node with the smallest such sum will serve as the total center.  The total connection cost of the tree that needs to be minimized is given by
\begin{equation*}
\ph_1(x^1, \ldots, x^k):=\sum_{i=1}^m \min_{\ell =1, \ldots, k} \sigma_F(x^\ell-a^i)\;\;+\;\;\min_{i=1, \ldots, m}\sum_{\ell =1}^{k}\sigma_F(x^\ell - a^i).
\end{equation*}
On the other hand, in the second model the $m$ nodes are clustered around $k+1$ artificial centers by trying to minimize the minimum sum of the distances from each artificial center to the remaining $k$ centers. Such a center will eventually be named as the total center. In this case, the total connection cost of the tree that needs to be minimized is given by
\begin{equation*}
\ph_2(x^1, \ldots, x^{k+1}):=\sum_{i=1}^m \min_{\ell =1, \ldots, k+1} \sigma_F(x^\ell-a^i) \;\;+\;\;\min_{\ell=1, \ldots, k+1}\sum_{j=1}^{k+1}\sigma_F( x^\ell - x^j).
\end{equation*}
The main difference between Model~I and Model~II is the way in which  the total center is selected. In addition, in Model II the total center also serves as a cluster center.

%\begin{figure}[H]
%\begin{subfigure}{.50\textwidth}
%  \centering
%  \includegraphics[width=0.95\textwidth]{Model1.png}
%  \caption{Model 1}
%  \label{fig:DS18}
%\end{subfigure}%
%\begin{subfigure}{.50\textwidth}
%  \centering
%\includegraphics[width=0.95\textwidth]{Model2.png}
%  \caption{Model 2}
%  \label{fig:DS1002}
%\end{subfigure}
%\caption{Plots of the two different models of a multicast network.}
%\label{fig:models}
%\end{figure}

The algorithms we will develop are expected to solve the continuous optimization models in a reasonable amount of time and give us   approximate solutions to the original discrete optimization models. Note that each node $a^i$ is assigned to its closest center $x^\ell$, but in both models the centers might not be real nodes. Therefore, for the continuous optimization model to solve (or approximate) the discrete model, we need to add a constraint that tries to minimize the difference between the artificial centers and the real centers, i.e.,
\begin{equation*}
\phi_1(x^1, \ldots, x^k):=\sum_{\ell=1}^k \min_{i=1, \ldots,m}\sigma_F(x^\ell-a^i)=0
\end{equation*}
and
\begin{equation*}
\phi_2(x^1, \ldots, x^{k+1}):=\sum_{\ell=1}^{k+1} \min_{i=1, \ldots,m}\sigma_F(x^\ell-a^i)=0.
\end{equation*}
Note that we use the generalized distance generated by $\sigma_F$ in the constraints for convenience of presentation although it is possible to use different distances such as the Euclidean distance.

Model I was originally proposed in \cite{bjor03} where the authors used the \emph{derivative-free discrete gradient method} established in \cite{B99} to solve the resulting optimization problem, but this method is not suitable for large-scale settings in high dimensions. It is also considered in \cite{am07} to solve a similar model where the \emph{squared Euclidean distance} used as a similarity measure. Model II was considered in \cite{bvx13} without constraints, and the \emph{hyperbolic smoothing technique} was used to solve the problem.

\section{Basic Definitions and Tools of Optimization}
\setcounter{equation}{0}

In this section, we  present two main tools of optimization used to solve the bilevel hierarchical crusting problem: the DCA introduced by Pham Dinh Tao and the Nesterov smoothing technique.

We consider throughout the paper DC programming:
\begin{equation}\label{dcf}
\mbox{\rm minimize}\, f(x):=g(x)-h(x), x\in \R^n,
\end{equation}
where $g\colon\R^n\to \R$ and $h\colon \R^n\to \R$ are convex functions. The function $f$ in \eqref{dcf} is called a \emph{DC function} and $g-h$ is called a \emph{DC decomposition} of $f$.

Given a convex function $g\colon \R^n\to \R$, the \emph{Fenchel conjugate} of $g$ is defined by
\begin{equation*}
g^*(y):=\sup\{\la y, x\ra -g(x)\; |\; x\in \R^n\}.
\end{equation*}
Note that  $g^*\colon \R^n\to (-\infty, +\infty]$ is also a convex function. In addition,  $x\in \partial g^*(y)$ if and only if $y\in \partial g(x)$, where $\partial$ denotes the subdifferential operator in the sense of convex analysis; see, e.g., \cite{hl93,BMN14B,r70b}.

Let us present below the DCA introduced by Tao and An \cite{at97,ta98} as applied to (\ref{dcf}). Although the algorithm is used for nonconvex optimization problems, the convexity of the functions involved still plays a crucial role.
\begin{algorithm}
  \caption{{\bf The DCA}}
  \begin{algorithmic}[1]
    \State \textbf{Input:} $x_0\in \R^n$, $N\in \N$.
    \State {\bf for} $k=1, \ldots, N$ {\bf do}
    \State \qquad Find $y_k\in \partial h(x_{k-1})$
    \State \qquad Find $x_{k}\in \partial g^*(y_k)$
    \State {\bf end for}
    \State \textbf{Output:} $x_N$.
  \end{algorithmic}
\end{algorithm}

Let us discuss below a convergence result of DC programming. A  function $h\colon \R^n\to \R$ is called $\gamma$-convex ($\gamma\geq 0$) if the function defined by $k(x):=h(x)-\frac{\gamma}{2}\|x\|^2$, $x\in \R^n$, is convex. If there exists $\gamma>0$ such that $h$ is $\gamma-$convex, then $h$ is called strongly convex.  We say that an element $\ox\in \R^n$ is a \emph{critical point} of the function $f$ defined by \eqref{dcf} if
\begin{equation*}\partial g(\ox)\cap \partial h(\ox)\neq \emptyset.
\end{equation*}
Obviously, in the case where both $g$ and $h$ are differentiable, $\ox$ is a critical point of $f$ if and only if $\ox$ satisfies the Fermat rule $\nabla f(\ox)=0$. The theorem below provides a convergence result for the DCA. It can be derived directly from \cite[Theorem 3.7]{ta98}.

\begin{Theorem} Consider the function $f$ defined by \eqref{dcf} and the sequence $\{x_k\}$ generated by the Algorithm 1. Then the following properties are valid:\\[1ex]
{\rm\bf(i)} If $g$ is $\gamma_1$-convex and $h$ is $\gamma_2$-convex, then
\begin{equation*}\label{sc1}
f(x_{k})-f(x_{k+1})\geq \frac{\gamma_1+\gamma_2}{2}\|x_{k+1}-x_k\|^2\; \mbox{\rm for all }k\in \N.
\end{equation*}
{\rm\bf (ii)} The sequence $\{f(x_k)\}$ is monotone decreasing.\\[1ex]
 {\rm\bf (iii)} If $f$ is bounded from below,   $g$ is $\gamma_1$-convex and $h$ is $\gamma_2$-convex with $\gamma_1+\gamma_2>0$, and  $\{x_k\}$ is bounded, then every subsequential limit of the sequence $\{x_k\}$ is a critical point of $f$.
\end{Theorem}

Let us present below a direct consequence of the Nesterov smoothing technique given in \cite{ny05}. In the proposition below, $d(x; \Omega)$ denotes the Euclidean distance and $P(x; \Omega)$ denotes the Euclidean projection from a point $x$ to a nonempty closed convex set $\Omega$ in $\R^n$.
\begin{Proposition}\label{p1} Given any $a\in \R^n$ and $\mu>0$, a Nesterov smoothing approximation of $\ph(x):=\sigma_F(x-a)$ has the representation
\begin{equation*}\label{NA}
\ph_\mu(x):=\frac{1}{2\mu}\|x-a\|^2-\frac{\mu}{2}\big[d(\frac{x-a}{\mu}; F)\big]^2.
\end{equation*}
Moreover, $\nabla \ph_\mu(x)=P(\frac{x-a}{\mu}; F)$ and
\begin{equation*}\label{approx}
\ph_\mu(x)\leq \ph(x)\leq \ph_\mu(x)+\frac{\mu}{2}\|F\|^2,
\end{equation*}
where $\|F\|:=\sup\{\|f\|\; |\; f\in F\}$.
\end{Proposition}

\section{Hierarchical Clustering via Continuous Optimization Techniques: Model I}\label{S:HC}
\setcounter{equation}{0}

In this section, we present an approach of using continuous optimization techniques for hierarchical clustering. As mentioned earlier, our main tools are the DCA and the Nesterov smoothing technique. Recall that the first model under consideration is formulated as a constrained optimization problem:

\begin{equation*}
\begin{aligned}
& \text{minimize}
& & \sum_{i=1}^m \min_{\ell =1, \ldots, k} \sigma_F(x^\ell-a^i)+\min_{i=1, \ldots, m}\sum_{\ell =1}^{k}\sigma_F(x^\ell - a^i)\\
& \text{subject to}
& & \sum_{\ell=1}^k \min_{i=1, \ldots,m}\sigma_F(x^\ell-a^i)=0,\; x^1, \ldots, x^k \in \mathbb{R}^n.\\
\end{aligned}
\end{equation*}
After the centers $x^1, \ldots, x^k$ have been found, a total center is selected from the existing nodes as follows: For each $i=1, \ldots, m$, we compute the sum $\sum_{\ell=1}^{k}\sigma_F(x^\ell - a^i)$. Then a total center $c^*$ is a node $a^i$ that yields the smallest sum, i.e.,
\begin{equation*}
c^*:= \mbox{\rm argmin}\Big\{\sum_{\ell=1}^{k}\sigma_F(x^\ell - a^i)\; \big| \;i=1, \ldots, m\Big\}.
\end{equation*}

Now we convert the constrained optimization problem under consideration to an unconstrained optimization problem using the penalty method with a penalty parameter $\lambda>0$:
\begin{align*}
  \mbox{\rm minimize }&\sum_{i=1}^{m}\min_{\ell=1, \ldots, k}\sigma_F(x^\ell - a^i) +  \min_{i=1, \ldots, m}\sum_{\ell=1}^{k}\sigma_F(x^\ell - a^i)+ \lambda \sum_{\ell=1}^{k}\min_{i=1, \ldots, m}\sigma_F(x^\ell - a^i)\\
  & x^1, \ldots, x^k \in \mathbb{R}^n.
  \end{align*}

\begin{Proposition}\label{p2} The objective function
\begin{equation*}
 f(x^1, \ldots, x^k):=\sum_{i=1}^{m}\min_{\ell=1, \ldots, k}\sigma_F(x^\ell - a^i) +  \min_{i=1, \ldots, m}\sum_{\ell=1}^{k}\sigma_F(x^\ell - a^i)+ \lambda \sum_{\ell=1}^{k}\min_{i=1, \ldots, m}\sigma_F(x^\ell - a^i)
 \end{equation*}
 for $x^1, \ldots, x^k \in \mathbb{R}^n$ and $\lambda>0$ can be written as a difference of convex functions.
\end{Proposition}
{\bf Proof.}  First note that the minimum of $m$ real numbers $\alpha_i$ for $i=1, \ldots, m$ has the representation:
\begin{equation*}
\min_{i=1,\ldots ,m}\alpha_i=\sum_{i=1}^{m} \alpha_i  -  \max_{t=1,\ldots ,m}\sum_{\substack{i=1\\ i\neq t}}^{m} \alpha_i.
\end{equation*}
Hence, we can represent $ f(x^1, \ldots, x^k)$ as a function defined on $(\R^n)^k$ as follows:
\begin{align*}
 f(x^1, \ldots, x^k)&=(2+\lambda) \sum_{i=1}^{m}\sum_{\ell=1}^{k}\sigma_F(x^\ell - a^i) -\sum_{i=1}^{m}\max_{t=1, \ldots, k}\sum_{\substack{\ell=1\\ \ell\neq t}}^{k}\sigma_F(x^\ell - a^i)\\
&-\lambda\sum_{\ell=1}^{k}\max_{t=1, \ldots, m}\sum_{\substack{i=1\\ i\neq t}}^{m}\sigma_F(x^\ell - a^i)- \max_{t=1, \ldots, m}\sum_{\substack{i=1\\ i\neq t}}^{m}\sum_{\ell=1}^{k}\sigma_F(x^\ell - a^i).
\end{align*}
This shows that $f$ has a DC  representation $f=g_0-h_0$, where
\begin{equation}\label{g0}
g_0(x^1, \ldots, x^k):=(2+\lambda) \sum_{i=1}^{m}\sum_{\ell=1}^{k}\sigma_F(x^\ell - a^i)
\end{equation}
and
\begin{align*}
h_0(x^1, \ldots, x^k)&:=\sum_{i=1}^{m}\max_{t=1, \ldots, k}\sum_{\substack{\ell=1\\ \ell\neq t}}^{k}\sigma_F(x^\ell - a^i)+\lambda\sum_{\ell=1}^{k}\max_{t=1, \ldots, m}\sum_{\substack{i=1\\ i\neq t}}^{m}\sigma_F(x^\ell - a^i)\\
&+ \max_{t=1, \ldots, m}\sum_{\substack{i=1\\ i\neq t}}^{m}\sum_{\ell=1}^{k}\sigma_F(x^\ell - a^i)
\end{align*}
are convex functions defined on $(\R^n)^k$. $\h$

Based on Proposition \ref{p1}, we obtain a Nesterov's approximation of the generalized  distance function $\ph(x):=\sigma_F(x-a)$ for $x, a\in \R^n$ as follows
\begin{equation*}
\ph_{\mu}(x): = \frac{\mu}{2}\left[\Bigg\|\frac{x-a}{\mu}\Bigg\|^2-\left[d\left(\frac{x- a}{\mu}; F\right)\right]^2\right].
\end{equation*}
As a result, the function $g_0$ defined in \eqref{g0} has a smooth approximation given by
\begin{equation*}
g_{0\mu}(x^1, \ldots, x^k):=\frac{(2+\lambda)\mu}{2} \sum_{i=1}^{m}\sum_{\ell=1}^{k}\Bigg\|\frac{x^\ell - a^i}{\mu}\Bigg\|^2 \;-\; \frac{(2+\lambda )\mu}{2} \sum_{i=1}^{m}\sum_{\ell=1}^{k}\left[d\left(\frac{x^\ell - a^i}{\mu}; F\right)\right]^2.
\end{equation*}
Thus, the function $f$ has the following DC approximation convenient for applying the DCA:
\begin{align*}
f_{\mu}(x^1, \ldots, x^k): &=\frac{(2+\lambda)\mu}{2} \sum_{i=1}^{m}\sum_{\ell=1}^{k}\Bigg\|\frac{x^\ell - a^i}{\mu}\Bigg\|^2 - \frac{(2+\lambda )\mu}{2} \sum_{i=1}^{m}\sum_{\ell=1}^{k}\left[d\left(\frac{x^\ell - a^i}{\mu}; F\right)\right]^2\\
&-\;\sum_{i=1}^{m}\max_{t=1, \ldots, k}\sum_{\substack{\ell=1\\ \ell\neq t}}^{k}\sigma_F(x^\ell - a^i)-  \lambda\;\sum_{\ell=1}^{k}\max_{s=1, \ldots, m}\sum_{\substack{i=1\\ i\neq t}}^{m}\sigma_F(x^\ell - a^i)\\
&- \max_{t=1, \ldots, m}\sum_{\substack{i=1\\ i\neq t}}^{m}\sum_{\ell=1}^{k}\sigma_F(x^\ell - a^i).
\end{align*}
Instead of minimizing the function $f$, we minimize its DC approximation
\begin{equation*}
 f_{\mu}(x^1, \ldots, x^k)=g_{\mu}(x^1, \ldots, x^k) - h_{\mu}(x^1, \ldots, x^k),\;\;x^1, \ldots, x^k \in \mathbb{R}^n.
 \end{equation*}
In this formulation, $g_{\mu}$ and $h_{\mu}$ are convex functions  given by
\begin{align*}
& g_{\mu}(x^1, \ldots, x^k):= \frac{2+\lambda}{2\mu} \sum_{i=1}^{m}\sum_{\ell=1}^{k}\|x^\ell - a^i\|^2,\\
& h_{\mu}(x^1, \ldots, x^k):= h_{1\mu}(x^1, \ldots, x^k) + h_{2}(x^1, \ldots, x^k) + h_{3}(x^1, \ldots, x^k) + h_{4}(x^1, \ldots, x^k),
\end{align*}
where
\begin{align*}
h_{1\mu}(x^1, \ldots, x^k) &:= \frac{(2+\lambda)\mu}{2} \sum_{i=1}^{m}\sum_{\ell=1}^{k}\left[d\left(\frac{x^\ell - a^i}{\mu}; F\right)\right]^2, \;  h_{2}(x^1, \ldots, x^k):= \sum_{i=1}^{m}\max_{t=1, \ldots, k}\sum_{\substack{\ell=1\\ \ell\neq t}}^{k}\sigma_F(x^\ell - a^i),\\
h_{3}(x^1, \ldots, x^k) &:= \lambda\;\sum_{\ell=1}^{k}\max_{t=1, \ldots, m}\sum_{\substack{i=1\\ i\neq t}}^{m}\sigma_F(x^\ell - a^i),\; h_{4}(x^1, \ldots, x^k): = \;\max_{t=1, \ldots, m}\sum_{\substack{i=1\\ i\neq t}}^{m}\sum_{\ell=1}^{k}\sigma_F(x^\ell - a^i).
\end{align*}
The proposition below is a direct consequence of Proposition \ref{p1}.
\begin{Proposition}\label{pp} Given any $\lambda > 0$ and $\mu > 0$, the functions $f$ and $f_{\mu}$ satisfy
\begin{equation*}\label{approx}
f_{\mu}(x^1, \ldots, x^k)\leq f(x^1, \ldots, x^k)\leq f_{\mu}(x^1, \ldots, x^k)+mk\left(1+\frac{\;\lambda\;}{2}\right)\mu\|F\|^2.
\end{equation*}
for all $x^1, \ldots, x^k\in \R^n$.
\end{Proposition}

In what follows we will prove that each of the functions $f$ and $f_\mu$ admits an absolute minimum in $(\R^n)^k$.
\begin{Theorem} Given any $\lambda > 0$ and $\mu > 0$, each of the functions $f$ and $f_{\mu}$ has an absolute minimum in $(\R^n)^k$.
\end{Theorem}
{\bf Proof.} Let us show that for any $\gamma\in \R$, the sublevel set
\begin{equation*}
\mathcal{L}_\gamma:=\{(x^1, \ldots, x^k)\; |\; f(x^1, \ldots, x^k)\leq \gamma\}
\end{equation*}
is bounded in $(\R^n)^k$. Since $0\in \mbox{\rm int}(F)$, there exists $r>0$ such that $\B(0; r)\subset F$. Consequently,
\begin{equation*}
r\|x\|=\sup\{\la x, u\ra\; |\; u\in \B(0; r)\}\leq \sup\{\la x, u\ra\; |\; u\in F\}=\sigma_F(x)\; \mbox{\rm for all }x\in \R^n.
\end{equation*}
From the definition of the function $f$, we have
\begin{align*}
\{(x^1, \ldots, x^k)\in (\R^n)^k\; |\; f(x^1, \ldots, x^k)\leq \gamma\}&\subset \{(x^1, \ldots, x^k)\in (\R^n)^k\; |\; \min_{i=1, \ldots, m}\sum_{\ell=1}^{k}\sigma_F(x^\ell - a^i)\leq \gamma\}\\
&\subset \{(x^1, \ldots, x^k)\in (\R^n)^k\; |\; \min_{i=1, \ldots, m}\sum_{\ell=1}^{k}\|x^\ell - a^i\|\leq \frac{\gamma}{r}\}\\
&\subset\bigcup_{i=1}^m\{(x^1, \ldots, x^k)\; |\; \ph_i(x^1, \ldots, x^k)\leq \frac{\gamma}{r}\},
\end{align*}
where $\ph_i(x^1, \ldots, x^k):=\sum_{\ell=1}^{k}\|x^\ell - a^i\|$. Observe that for each $i=1, \ldots, m$, one has the inclusion
\begin{equation*}
\{(x^1, \ldots, x^k)\; |\; \ph_i(x^1, \ldots, x^k)\leq \frac{\gamma}{r}\}\subset \{(x^1, \ldots, x^k)\; |\; \sum_{\ell=1}^k\|x^\ell\|\leq \frac{\gamma}{r}+k\|a^i\|\}.
\end{equation*}
Thus, $\mathcal{L}_\gamma$ is a bounded set as it is contained in the union of a finite number of bounded sets in $(\R^n)^k$. As $f$ is a continuous function, it has an absolute minimum in $(R^n)^k$.

Let $\gamma_\mu:=mk\left(1+\frac{\;\lambda\;}{2}\right)\mu\|F\|^2$. It follows from Proposition \ref{pp} that for any $\gamma\in \R$,
\begin{equation*}
\{(x^1, \ldots, x^k)\in (\R^n)^k\; |\; f_\mu(x^1, \ldots, x^k)\leq \gamma\}\subset \{(x^1, \ldots, x^k)\in (\R^n)^k\; |\; f(x^1, \ldots, x^k)\leq \gamma_\mu+\gamma\}.
\end{equation*}
It follows that the sublevel set $\{(x^1, \ldots, x^k)\in (\R^n)^k\; |\; f_\mu(x^1, \ldots, x^k)\leq \gamma\}$ is also bounded, and hence $f_\mu$ has an absolute minimum in $(\R^n)^k$. $\h$

To facilitate the gradient and subgradient calculations for the DCA, we will introduce a \emph{data matrix} $\mathbf A$ and a \emph{variable matrix} $\mathbf X$.  The data matrix $\mathbf A$ is formed by putting each $a^i$, $i =1, \ldots, m$, in the $i^{th}$ row, i.e.,
\begin{equation*}
\mathbf{A} =
   \left(
   \begin{matrix}
      a_{11}      	&    a_{12}   	&    a_{13}		& \dots   		& a_{1n}\\
      a_{21}      	&    a_{22}   	&    a_{23}		& \dots   		& a_{2n}\\
      \vdots 		&  \vdots 		& \vdots  		& 		 	& \vdots\\
      a_{m1}      	&    a_{m2}   	&    a_{m3}	& \dots   		& a_{mn}
   \end{matrix}
   \right).
\end{equation*}
Similarly, if $x^1, \ldots, x^k$ are the $k$ cluster centers, then the variable $\mathbf X$ is formed by putting each $x^\ell$, $\ell =1, \ldots ,k$, in the $\ell^{th}$ row, i.e.,
\begin{equation*}
\mathbf{X} =
   \left(
   \begin{matrix}
      x_{11}      &    x_{12}   &    x_{13}& \dots   & x_{1n}\\
      x_{21}      &    x_{22}   &    x_{23}& \dots   & x_{2n}\\
      \vdots 		&  \vdots 		& \vdots  	& 	& \vdots\\
      x_{k1}      &    x_{k2}   &    x_{k3}& \dots   & x_{kn}
   \end{matrix}
   \right).
\end{equation*}
With these notations, the decision variable $\mathbf X$ of the optimization problem belongs to $\mathbb{R}^{k\times n}$, the linear space of $k\times n$ real matrices. Hence, we will assume that $\mathbb{R}^{k\times n}$ is equipped with the inner product $\langle X, Y\rangle:=\mbox{\rm trace}(X^TY)$. The \emph{Frobenius norm} on $\mathbb{R}^{k\times n}$ is defined by
\begin{equation*}
\|\mathbf X\|_F := \sqrt{\langle \mathbf X, \mathbf X \rangle} = \sqrt{ \sum_{\ell =1}^{k}\langle x^\ell, x^\ell\rangle} = \sqrt{\sum_{\ell =1}^{k}\|x^\ell\|^2}.
\end{equation*}
Let us start by computing the gradient of the first  part  of the DC decomposition, i.e.,
\begin{equation*}
g_{\mu}(\mathbf X) \;=\; \frac{2+\lambda}{2\mu} \sum_{i=1}^{m}\sum_{\ell=1}^{k}\|x^\ell - a^i\|^2.
\end{equation*}
Using the Frobenius norm, the function $g_{\mu}$  can  be written as
\begin{align*}
g_{\mu}(\mathbf X) &=\; \frac{2+\lambda}{2\mu} \sum_{i=1}^{m}\sum_{\ell =1}^{k}\|x^\ell - a^i\|^2\\
					&=\; \frac{2+\lambda}{2\mu} \sum_{i=1}^{m}\sum_{\ell =1}^{k}\left[\|x^\ell\|^2-2\langle x^\ell,  a^i\rangle  + \|a^i\|^2\right]\\
					&=\; \frac{2+\lambda}{2\mu} \left[m\|\mathbf X\|_F^2  -  2\langle \mathbf X, \mathbf E\mathbf A\rangle  +  k\|\mathbf A\|_F^2\right],
\end{align*}
where $\mathbf E$ is a $k\times m$ matrix whose entries are all ones. Hence,  $g_{\mu}$ is differentiable and its gradient is given by
\[\nabla g_{\mu}(\mathbf X)= \frac{2+\lambda}{\mu}\left[m\mathbf X - \mathbf E\mathbf A\right].\]															
Our goal now is to find  $\mathbf X\in\partial g^*(\mathbf Y)$, which can be accomplished by employing the relation
\begin{equation*}
\mathbf X\in\partial g^*(\mathbf Y)\; \mbox{\rm if and only if }\mathbf Y\in\partial g(\mathbf X).
\end{equation*}
This can equivalently be written as $\frac{2+\lambda}{\mu}\left[m\mathbf X - \mathbf E\mathbf A\right] =\mathbf Y$, and we solve for $\mathbf X$ as follows:
\begin{align*}
&(2+\lambda)\left[m\mathbf X - \mathbf E\mathbf A\right] =\mu \mathbf Y\\
&(2+\lambda)\mathbf X=  (2+\lambda)\mathbf E \mathbf A +\mu \mathbf Y\\
&\mathbf X 	= \frac{(2+\lambda)\mathbf E \mathbf A +\mu \mathbf Y}{\;\;(2+\lambda)m}
\end{align*}

Next, we will demonstrate in more detail the techniques we used to compute a subgradient for the convex function
\begin{equation*}
h_{\mu} = h_{1\mu}+\sum_{j=2}^{4}h_{j}.
\end{equation*}
Since each function in this sum is convex, we will compute a subgradient of $h_\mu$ applying the subdifferential sum rule (see, e.g.,  \cite[Corollary 2.46]{BMN14B}) and maximum rule (see, e.g.,  \cite[Proposition 2.54]{BMN14B}) well known in convex analysis. We will begin our demonstration with $h_{1\mu}$ given by
\begin{equation*}
h_{1\mu}(\mathbf X)  = \frac{(2+\lambda)\mu}{2} \sum_{i=1}^{m}\sum_{\ell=1}^{k}\left[d\left(\frac{x^\ell - a^i}{\mu}; F\right)\right]^2.
\end{equation*}
From its representation one can see that $h_{1\mu}$ is differentiable. Thus, its gradient at $\mathbf X$ can be computed by computing the \emph{partial derivatives} with respect to $x^1, \ldots, x^k$, i.e.,
\begin{equation}\label{partial derivative}
\frac{\partial h_{1\mu}}{\partial x^\ell }(\mathbf X)  = (2+\lambda) \sum_{i=1}^{m}\left[\frac{x^\ell - a^i}{\mu}\; - P\left(\frac{x^\ell - a^i}{\mu}; F \right)\right]\;\; \text{for} \;\;\ell = 1, \ldots, k.
\end{equation}
Hence, $\nabla h_{1\mu}(\mathbf X))$ is a $k\times n$ matrix $\mathbf H_{1}$ whose $\ell^{th}$ row is $\frac{\partial h_{1\mu}}{\partial x^\ell }(\mathbf X))$.

Note that the convex functions $h_{j}$ for $j=2, 3, 4$ are not differentiable in general. However, we can compute a subgradient for each function at $\mathbf X$ by applying the subdifferential sum rule and maximum rule for convex functions. The following is an illustration of how one can compute  subgradients of such functions using $h_{2}$ as an example. For $t=1, \ldots, k$ and $i=1, \ldots, m$, define
\begin{equation*}
\gamma_{ti}(\mathbf X):=\sum_{\ell=1, \ell\neq t}^{k}\sigma_F(x^\ell - a^i)=\sum_{\ell=1}^k\sigma_F(x^\ell - a^i)-\sigma_F(x^t-a^i)\;\text{and}\;\gamma_i(\mathbf X):=\max_{t=1,\ldots,k} \gamma_{ti}(\mathbf X).
\end{equation*} Thus, $h_{2}$ can be represented as the sum of $m$ convex functions as follows:
\begin{equation*}
h_{2}(\mathbf X) = \sum_{i=1}^{m}\max_{t=1, \ldots, k}\sum_{\ell=1, \ell\neq t}^{k}\sigma_F(x^\ell - a^i)=\sum_{i=1}^m\gamma_i(\mathbf X).
\end{equation*}
Note that $\gamma_i$ is the maximum of $k$ convex functions $\gamma_{ti}$ for $t=1, \ldots, k$. Based on the subdifferential maximum rule, for each $i=1, \ldots, m$, we will find
a $k\times n$ matrix $\mathbf H_{2i}\in \partial \gamma_i(\mathbf X)$. Then, by the subdifferential sum rule  $\mathbf H_{2}:=\sum_{i=1}^m \mathbf H_{2i}$ is a  subgradient of  $h_{2}$ at $\mathbf X$. To accomplish this goal, we first choose an index $t^*\in\{1, \ldots, k\}$ such that $\gamma_i(\mathbf X)=\gamma_{t^*i}(\mathbf X):=\sum_{\ell=1, \ell\neq t^*}^{k}\sigma_F(x^\ell - a^i)$. The $\ell^{\rm th}$ row $w_\ell^i$ of the matrix $\mathbf H_{2i}$ for $\ell\neq t^*$ can be computed as described in Proposition \ref{p2} below, which follows from \cite[Theorem 2.93]{BMN14B}. The $t^*$ row of the matrix $\mathbf H_{2i}$ is set to zero, as $\gamma_{it^*}$ is independent of $x^{t^*}$. The procedures for computing a subgradient for $h_{3}$ and $h_{4} $ are very similar to the procedure we have illustrated.

\begin{Proposition}\label{p2} Given $a\in \R^n$, the function $\ph(x):=\sigma_F(x-a)$ is convex with its subdifferential at $\ox\in \R^n$ given by
\begin{equation*}
\partial \ph(\ox)=\mbox{\rm co}\,F(\ox),
\end{equation*}
where $F(\ox):=\{q\in F\; |\; \la \ox, q\ra=\sigma_F(\ox)\}$.

In particular, if $F$ is the Euclidean closed unit ball in $\R^n$, then
%\begin{equation*}
\[
\partial \ph(\ox)=\Bigg\{\begin{array}{cr}
 \frac{\ox-a}{\|\ox-a\|}        &\text{if}\;\ox\neq a,\\
 \mathbb B        &\text{if}\;\ox=a.
\end{array}
\]
%\end{equation*}
\end{Proposition}

At this point, we have demonstrated all the necessary steps in calculating the gradients and subgradients needed for our first DCA-based algorithm for solving the bilevel hierarchical clustering problem formulated in Model~I.

\begin{algorithm}
  \caption{$\;$ {\bf Model~I}}
  \begin{algorithmic}[1]
    \State \textbf{Input:} $ \mathbf A, \mathbf X_0, \lambda_0, \mu_0,\sigma_1, \sigma_2, \epsilon, N \in \mathbb{N}$.
    \State {\bf while} \mbox{stopping criteria} ($\lambda$, $\mu$, $\epsilon$) = false {\bf do}
    \State \qquad {\bf for} $k = 1,\ldots, N\;\;$ {\bf do}
    \State \qquad{\qquad Find $\mathbf Y_k\in \partial h_{\mu}(\mathbf X_{k-1})$}
    \State \qquad{\qquad $\mathbf X_k = \frac{(2+\lambda)\mathbf E\mathbf A +\mu \mathbf Y_k}{\;\;(2+\lambda)m\;\;}$}
    \State \qquad {\bf end for}
    \State \qquad update $\lambda \;\;\text{and} \; \;\mu$
    \State {\bf end while}
    \State \textbf{Output:} $x_N$.
  \end{algorithmic}
\end{algorithm}

\begin{Example}\label{ex1}{\rm {\bf\rm ($\ell^2-$clustering with Algorithm 2).} In this example, we illustrate our method to study the problem of $\ell^2-$clustering. The key point in Algorithm 2 is the computation of $\mathbf Y\in \partial h_{\mu}(\mathbf X)$ for the case where $F$ is the Euclidean closed unit ball $\B$ in $\R^n$. By the subdifferential sum rule,
\begin{equation*}
h_\mu(\mathbf X)=\nabla h_{1\mu}(\mathbf X)+\partial h_2(\mathbf X)+\partial h_3(\mathbf X)+\partial h_4(\mathbf X).
\end{equation*}
Define
\begin{equation*}
u_{\ell i}:=\begin{cases}
      \frac{x^\ell-a^i}{\|x^\ell-a^i\|} &\mbox{\rm if}\;  x^\ell\neq a^i, \\
            0 & \mbox{\rm otherwise}.
   \end{cases}
   \end{equation*}
  Now, we illustrate the way to find the gradient of $h_1$ and a subgradient of $h_{i}$ for $i=2, 3, 4$ at $\mathbf X$.\\[1ex]
  {\em The gradient of $h_{1}$}: The gradient $\mathbf Y_1:=\nabla h_1(\mathbf X)$ is the $k\times n$ matrix whose $\ell^{th}$ row is $\frac{\partial h_{1\mu}}{\partial x^\ell }(\mathbf X)$ given in \eqref{partial derivative}. Note that in this case, the Euclidean projection $P(z; F)$ from $z\in \R^n$ to $F$ is given by
  \begin{equation*}
P(z; F):=\begin{cases}
      \frac{z}{\|z\|} &\mbox{\rm if}\;  \|z\|>1, \\
            z & \mbox{\rm otherwise}.
   \end{cases}
   \end{equation*}
  {\em A subgradient of $h_{2}$}: In this case,
  \begin{equation*}
  h_{2}(\mathbf X)=\sum_{i=1}^{m}\max_{t=1, \ldots, k}\sum_{\substack{\ell=1\\ \ell\neq t}}^{k}\|x^\ell - a^i\|=\sum_{i=1}^{m}\max_{t=1, \ldots, k}\big(\sum_{\ell=1}^k\|x^\ell-a^i\|-\|x^t-a^i\|\big).
  \end{equation*}
  For each $i=1, \ldots, m$, choose an index $t(i)$ such that
  \begin{equation*}
  \max_{t=1, \ldots, k}\big(\sum_{\ell=1}^k\|x^\ell-a^i\|-\|x^t-a^i\|\big)=\sum_{\ell=1}^k\|x^\ell-a^i\|-\|x^{t(i)}-a^i\|.
  \end{equation*}
  Let us now form a $k\times mn$ \emph{block matrix} $\mathbf U=(u_{\ell i})$, where $u_{\ell i}$ is considered  as a row vector.  We also use  $U^i$ to denote the $i^{\rm th}$ block column of the  matrix $\mathbf U$. Equivalently, $U_i$ is the $k\times n$ matrix formed by placing the row vectors $u_{\ell i}$ in its $\ell^{\rm th}$ row for $\ell=1, \ldots, k$.  Then a subgradient of $h_{2}$ at $\mathbf X$ is given by
  \begin{equation*}
  \mathbf Y_2:=\sum_{i=1}^m\big(U^i-e_{t(i)}u_{t(i) i}\big),
  \end{equation*}
  where $e_{t(i)}$ is the column vector of $k$ components with $1$ at the $t(i)^{\rm th}$ position and $0$ at other positions. \\[1ex]
  {\em A subgradient of $h_{3}$}: In this case,
  \begin{equation*}
 h_{3}(\mathbf X)= \lambda\;\sum_{\ell=1}^{k}\max_{t=1, \ldots, m}\sum_{\substack{i=1\\ i\neq t}}^{m}\|x^\ell - a^i\|=\lambda\;\sum_{\ell=1}^{k}\max_{t=1, \ldots, m}\big(\sum_{i=1}^m\|x^\ell - a^i\|-\|x^\ell-a^t\|\big).
  \end{equation*}
  For each $\ell=1, \ldots, k$, we choose an index $t(\ell)$ such that
  \begin{equation*}
  \max_{t=1, \ldots, m}\big(\sum_{i=1}^m\|x^\ell - a^i\|-\|x^\ell-a^t\|\big)=\sum_{i=1}^m\|x^\ell - a^i\|-\|x^\ell-a^{t(\ell)}\|.
  \end{equation*}
  Let $\mathbf V$ be the $k\times n$ matrix whose $\ell^{\rm th}$ row is $\sum_{i=1}^mu_{\ell i}-u_{\ell t(\ell)}$. Then a subgradient of $h_{3}$ at $X$ is given by
  \begin{equation*}
  \mathbf Y_3:=\lambda \mathbf V.
  \end{equation*}
   {\em A subgradient of $h_{4}$}: In this case,
   \begin{align*}
   h_{4}(\mathbf X)&=\max_{t=1, \ldots, m}\sum_{\substack{i=1\\ i\neq t}}^{m}\sum_{\ell=1}^{k}\|x^\ell - a^i\|\\
   &=\max_{t=1, \ldots, m}\left( \sum_{i=1}^m\sum_{\ell=1}^k\|x^\ell-a^i\|-\sum_{\ell=1}^k\|x^\ell-a^t\|\right).
   \end{align*}
   Again, we choose an index $t$ such that
   \begin{equation*}
   \max_{t=1, \ldots, m}\left( \sum_{i=1}^m\sum_{\ell=1}^k\|x^\ell-a^i\|-\sum_{\ell=1}^k\|x^\ell-a^t\|\right)=\sum_{i=1}^m\sum_{\ell=1}^k\|x^\ell-a^i\|-\sum_{\ell=1}^k\|x^\ell-a^t\|.
   \end{equation*}
   Let $\mathbf Z$ be the $k\times n$ matrix whose $\ell^{th}$ row is $\sum_{i=1}^m u_{\ell i}$. Then a subgradient of $h_{4}$ is given by
   \begin{equation*}
   \mathbf Y_4:=Z-Z_t,
   \end{equation*}
   where $Z_t$ is the $k\times n$ matrix whose $\ell^{\rm th}$ row is $u_{\ell t}$.

 } \end{Example}

\begin{Example}\label{ex2}{\rm {\bf\rm ($\ell^1-$clustering with Algorithm 2).} In this example, we illustrate our method to study the problem of $\ell^1-$clustering. We will find a subgradient $\mathbf Y\in \partial h_{\mu}(\mathbf X)$ for the case where $F$ is the \emph{closed unit box} in $\R^n$ given by
\begin{equation*}
F:=\{(u_1, \ldots, u_n)\in \R^n\; |\; -1\leq u_i\leq 1\; \mbox{\rm for }i=1, \ldots, n\}.
\end{equation*}
For $t\in \R$, define
\begin{equation*}
\mbox{\rm sign}(t):=\begin{cases}
      1 & t> 0, \\
      0 & t=0, \\
      -1 & t<0.
   \end{cases}
   \end{equation*}
  Then we define $\mbox{\rm sign}(x):=(\mbox{\rm sign}(x_1), \ldots, \mbox{\rm sign}(x_n))$ for $x=(x_1, \ldots, x_n)\in \R^n$. Note that for the function $p(x):=\|x\|_1$, a subgradient of $p$ at $x\in\R^n$ is simply $\mbox{\rm sign}(x)$. Now, we illustrate the way to find the gradient of $h_1$ and a subgradient of $h_{i}$ for $i=2, 3, 4$ at $\mathbf X$.\\[1ex]
   {\em The gradient of $h_{1}$}: Similar to Example \ref{ex1}, the gradient of $\mathbf Y_1:=\nabla h_1(\mathbf X)$ is the $k\times n$ matrix whose $\ell^{th}$ row is $\frac{\partial h_{1\mu}}{\partial x^\ell }(\mathbf X)$ given in \eqref{partial derivative}. Note that in this case, the Euclidean projection $P(z; F)$ from $z\in \R^n$ to $F$ is given by
  \begin{equation*}
P(z; F):=\max(-e, \min (z, e))\; \mbox{\rm componentwise},
   \end{equation*}
   where $e\in \R^n$ is the vector consisting of $1$ in each component. \\[1ex]
  {\em A subgradient of $h_{2}$}: In this case,
  \begin{equation*}
  h_{2}(\mathbf X)=\sum_{i=1}^{m}\max_{r=1, \ldots, k}\sum_{\substack{\ell=1\\ \ell\neq r}}^{k}\|x^\ell - a^i\|_1=\sum_{i=1}^{m}\max_{r=1, \ldots, k}\big(\sum_{\ell=1}^k\|x^\ell-a^i\|_1-\|x^r-a^i\|_1\big).
  \end{equation*}
  For each $i=1, \ldots, m$, choose an index $r(i)$ such that
  \begin{equation*}
  \max_{r=1, \ldots, k}\big(\sum_{\ell=1}^k\|x^\ell-a^i\|_1-\|x^r-a^i\|_1\big)=\sum_{\ell=1}^k\|x^\ell-a^i\|_1-\|x^{r(i)}-a^i\|_1.
  \end{equation*}
  Now we form the $k\times mn$ \emph{signed block matrix} $S=(s_{\ell i})$ given by $s_{\ell i}=\mbox{\rm sign}(x^\ell-a^i)$ as a row vector.  We also use $S^i$ to denote the $i$th column block matrix of the signed matrix $S$. Then a subgradient of $h_{2}$ at $\mathbf X$ is given by
  \begin{equation*}
  \mathbf Y_2:=\sum_{i=1}^m\big(S^i-e_{r(i)}s_{r(i) i}\big),
  \end{equation*}
  where $e_{r(i)}$ is the column vector of $k$ components with $1$ at the $r(i)$th position and $0$ at other positions. \\[1ex]
  {\em A subgradient of $h_{3}$}: In this case,
  \begin{equation*}
 h_{3}(\mathbf X)= \lambda\;\sum_{\ell=1}^{k}\max_{t=1, \ldots, m}\sum_{\substack{i=1\\ i\neq t}}^{m}\sigma_F(x^\ell - a^i)=\lambda\;\sum_{\ell=1}^{k}\max_{t=1, \ldots, m}\big(\sum_{i=1}^m\|x^\ell - a^i\|_1-\|x^\ell-a^t\|_1\big).
  \end{equation*}
  For each $\ell=1, \ldots, k$, we choose an index $t(\ell)$ such that
  \begin{equation*}
  \max_{t=1, \ldots, m}\big(\sum_{i=1}^m\|x^\ell - a^i\|_1-\|x^\ell-a^t\|_1\big)=\sum_{i=1}^m\|x^\ell - a^i\|_1-\|x^\ell-a^{t(\ell)}\|_1.
  \end{equation*}
  Let $\mathbf V$ be the $k\times n$ matrix whose $\ell^{\rm th}$   row is $\sum_{i=1}^ms_{\ell i}-s_{\ell t(\ell)}$. Then a subgradient of $h_{3}$ at $\mathbf X$ is given by
  \begin{equation*}
  \mathbf Y_3:=\lambda \mathbf V.
  \end{equation*}
   {\em A subgradient of $h_{4}$}: In this case,
   \begin{align*}
   h_{4}(\mathbf X)&=\max_{t=1, \ldots, m}\sum_{\substack{i=1\\ i\neq t}}^{m}\sum_{\ell=1}^{k}\|x^\ell - a^i\|_1\\
   &=\max_{t=1, \ldots, m}\left( \sum_{i=1}^m\sum_{\ell=1}^k\|x^\ell-a^i\|_1-\sum_{\ell=1}^k\|x^\ell-a^t\|_1\right).
   \end{align*}
   Again, we choose an index $t$ such that
   \begin{equation*}
   \max_{t=1, \ldots, m}\left( \sum_{i=1}^m\sum_{\ell=1}^k\|x^\ell-a^i\|_1-\sum_{\ell=1}^k\|x^\ell-a^t\|_1\right)=\sum_{i=1}^m\sum_{\ell=1}^k\|x^\ell-a^i\|_1-\sum_{\ell=1}^k\|x^\ell-a^t\|_1.
   \end{equation*}
   Let $T$ be the $k\times n$ matrix whose $\ell^{\rm th}$ row is $\sum_{i=1}^m s_{\ell i}$. Then a subgradient of $h_{4}$ is given by
   \begin{equation*}
   \mathbf Y_4:=T-T_t,
   \end{equation*}
   where $T_t$ is the $k\times n$ matrix whose $\ell^{\rm th}$ row is $s_{\ell t}$.

 } \end{Example}

\section{Hierarchical Clustering via Continuous Optimization Techniques: Model II}\label{S:HC}
\setcounter{equation}{0}

In this section, we focus on developing nonconvex optimization techniques based on the DCA and the Nesterov smoothing technique for the second model. Similar to Model I, we will solve the following constrained optimization problem:
\begin{equation*}
\begin{aligned}
& \text{minimize}
& & \sum_{i=1}^m \min_{\ell=1, \ldots , k+1} \sigma_F(x^\ell-a^i)+\min_{\ell=1, \ldots , k+1}\sum_{j =1}^{k+1}\sigma_F(x^\ell - x^j)\\
& \text{subject to}
& & \sum_{\ell=1}^{k+1} \min_{i=1, \ldots ,m}\sigma_F(x^\ell-a^i)=0,\; x^1, \ldots, x^{k+1} \in \mathbb{R}^n.\\
\end{aligned}
\end{equation*}
The total center is determined by
\begin{equation*}
c^*:= \mbox{\rm argmin }\Big\{\sum_{j=1}^{k+1}\sigma_F(x^\ell - x^j)\; |\; \ell=1, \ldots, k+1\Big\}.
\end{equation*}
This constrained optimization problem can be solved by the following unconstrained optimization problem by the penalty method with a penalty parameter $\lambda>0$:
\begin{align*}
 \text{minimize} & \sum_{i=1}^{m}\min_{\ell=1, \ldots,  k+1}\sigma_F(x^\ell - a^i) +  \min_{\ell=1, \ldots,  k+1}\sum_{j =1}^{k+1}\sigma_F(x^\ell - x^j)+ \lambda \sum_{\ell=1}^{k+1}\min_{i=1, \ldots, m}\sigma_F(x^\ell - a^i)\\
 & x^1, \ldots, x^{k+1} \in \mathbb{R}^n.
\end{align*}
With the Nesterov smoothing technique, the objective function has the following approximation that is convenient for implementing the DCA:
\begin{align*}
f_{\mu}(\mathbf X):&=\;\frac{(1+\lambda)\mu}{2} \sum_{i=1}^{m}\sum_{\ell=1}^{k+1}\Bigg\|\frac{x^\ell - a^i}{\mu}\Bigg\|^2 +\frac{\mu}{2} \sum_{\ell=1}^{k+1}\sum_{j=1}^{k+1}\Bigg\|\frac{x^\ell - x^j}{\mu}\Bigg\|^2\\
&-\; \frac{(1+\lambda)\mu}{2} \sum_{i=1}^{m}\sum_{\ell=1}^{k+1}\left[d\left(\frac{x^\ell - a^i}{\mu}; F\right)\right]^2\;-\;\sum_{i=1}^{m}\max_{r=1, \ldots, k+1}\sum_{\substack{\ell=1\\ \ell\neq r}}^{k+1}\sigma_F(x^\ell - a^i) \\
& -  \lambda\;\sum_{\ell=1}^{k+1}\max_{t=1, \ldots, m}\sum_{\substack{i=1\\ i\neq t}}^{m}\sigma_F(x^\ell - a^i)-\frac{\mu}{2} \sum_{\ell=1}^{k+1}\sum_{j=1}^{k+1}\left[d\left(\frac{x^\ell - x^j}{\mu}; F\right)\right]^2\;-\; \max_{r=1, \ldots, k+1}\sum_{\substack{\ell=1\\ \ell\neq r}}^{k+1}\sum_{j=1}^{k+1}\sigma_F(x^\ell - x^j).
\end{align*}
As in the previous section, we use a variable matrix $\mathbf X$ of size $(k+1)\times n$ to store the row vector $x^\ell$ in its $\ell^{\rm th}$ row for $\ell=1, \ldots, k+1$. Now we solve the following DC programming:
\begin{equation*}
\begin{aligned}
&\text{minimize}
&f_{\mu}(\mathbf X)=g_{\mu}(\mathbf X) - h_{\mu}(\mathbf X),\;\mathbf X\in \R^{(k+1)\times n},
\end{aligned}
\end{equation*} where $g_{\mu}$ and $h_{\mu}$ are convex functions by
\begin{equation}\label{gmu}
g_{\mu}(\mathbf X):= g_{1\mu}(\mathbf X) + g_{2\mu}(\mathbf X)
\end{equation}
 and
\[h_{\mu}(\mathbf X):= h_{1\mu}(\mathbf X) + h_{2\mu}(\mathbf X) + h_{3\mu}(\mathbf X) + h_{4\mu}(\mathbf X)+ h_{5\mu}(\mathbf X),\] where their respective components are defined as follows:
\begin{equation*}
 g_{1\mu}(\mathbf X):= \frac{(1+\lambda)\mu}{2} \sum_{i=1}^{m}\sum_{\ell=1}^{k+1}\Bigg\|\frac{x^\ell - a^i}{\mu}\Bigg\|^2,\; g_{2\mu}(\mathbf X) = \frac{\mu}{2} \sum_{\ell=1}^{k+1}\sum_{j=1}^{k+1}\Bigg\|\frac{x^\ell - x^j}{\mu}\Bigg\|^2
\end{equation*}
and
\begin{align*}
& h_{1\mu}(\mathbf X):=\frac{(1+\lambda)\mu}{2} \sum_{i=1}^{m}\sum_{\ell=1}^{k+1}\left[d\left(\frac{x^\ell - a^i}{\mu}; F\right)\right]^2, \;
h_{2\mu}(\mathbf X):  = \frac{\mu}{2} \sum_{\ell=1}^{k+1}\sum_{j=1}^{k+1}\left[d\left(\frac{x^\ell - x^j}{\mu}; F\right)\right]^2\\
& h_{3}(\mathbf X):= \sum_{i=1}^{m}\max_{t=1, \ldots, k+1}\sum_{\substack{\ell=1\\ \ell\neq t}}^{k+1}\sigma_F(x^\ell - a^i), \; h_{4}(\mathbf X):= \lambda\;\sum_{\ell=1}^{k+1}\max_{t=1, \ldots, m}\sum_{\substack{i=1\\ i\neq t}}^{m}\sigma_F(x^\ell - a^i), \\
& h_{5}(\mathbf X):= \max_{t=1, \ldots, k+1}\sum_{\substack{\ell=1\\ \ell\neq t}}^{k+1}\sum_{j=1}^{k+1}\sigma_F(x^\ell - x^j).
\end{align*}

\begin{Lemma}\label{lm51} Let $\mathbf E$ be square matrix with size $(k+1)$ whose entries are all ones  and let  $\mathbb{I}$ be the identity matrix of size $(k+1)$. \\[1ex]
{\rm\bf (i)} Given any  real numbers $a$ and $b$ with $a\neq 0$ and $a\neq -(k+1)b$, the matrix $\mathbf M:=a \Bbb I +b\mathbf E$ is invertible with
\begin{equation*}
\mathbf M^{-1}=x \mathbb{I}+y\mathbf E,
\end{equation*}
where  $x=\dfrac{1}{a}$ and $y=-\dfrac{b}{a[a+b(k+1)]}$.\\[1ex]
{\rm\bf (ii)} Let $\widetilde{\mathbf E}:= (k+1)\mathbb{I}-\mathbf E$.  Given any  real numbers $c$ and $d$ with $c\neq 0$ and $c\neq -d(k+1)$, the matrix $\mathbf N:=c\Bbb I+d\widetilde{\mathbf E}$ is invertible  with
\begin{equation*}
\mathbf N^{-1}=\alpha \mathbb{I}+\beta \mathbf E,
\end{equation*}
where $\alpha=\frac{1}{c+d(k+1)}$ and $\beta=\dfrac{d}{c[c+d(k+1)]}$.
\end{Lemma}
{\bf Proof.} {\rm\bf (i)} Observe that
\begin{align*}
(a \Bbb I +b\mathbf E)(x \Bbb I+y\mathbf E)&=ax\Bbb I+(bx+ay)\mathbf E+by\mathbf E^2\\
&=ax\Bbb I+(bx+ay)\mathbf E+by(k+1)\mathbf E.
\end{align*}
Thus, $(a \Bbb I +b\mathbf E)(x \Bbb I+y\mathbf E)=\Bbb I$ if and only if
\begin{equation*}
ax=1\; \mbox{\rm and }bx+[a+b(k+1)]y=0.
\end{equation*}
Equivalently, $x=\dfrac{1}{a}$ and $y=-\dfrac{b}{a[a+b(k+1)]}$.\\[1ex]
{\rm\bf (ii)} We have
\begin{equation*}
N=c\Bbb I+d\widetilde{\mathbf E}=[c+d(k+1)\Bbb I]-d\mathbf E.
\end{equation*}
It remains to apply the result from (i). $\h$

The proposition below provides a formula for computing $\nabla g^*_\mu$ required for applying the DCA.

\begin{Proposition} Given any $\lambda>0$ and $\mu>0$, the Fenchel conjugate $g_\mu^*$ of the function $g_\mu$ defined in \eqref{gmu} is continuously differentiable with
\begin{equation*}
\nabla g_\mu^*(\mathbf Y)=\left(\alpha\mathbb{I}+\beta \mathbf E\right)\Big((1+\lambda)\mathbf E \mathbf A + \mu \mathbf Y\Big)\; \mbox{\rm for }\mathbf Y\in \R^{k\times n},
\end{equation*}
where $\mathbf E$ is defined in Lemma \ref{lm51} and
\begin{equation}\label{ab}
\alpha:=\frac{1}{m(\lambda+1)+2(k+1)}\; \mbox{\rm and }\beta:=\frac{2}{m(\lambda+1)[m(\lambda+1)+2(k+1)]}.
\end{equation}
\end{Proposition}
{\bf Proof.}  We have
\begin{align*}
&\nabla g_{1\mu}(\mathbf X) \;=\;\frac{1+\lambda}{\mu}\left[m\mathbf X - \mathbf E\mathbf A\right],\\
&\nabla g_{2\mu}(\mathbf X) \;=\;\frac{2}{\mu}\left[(k+1)\mathbb {I}-\mathbf E\right]\mathbf X.
\end{align*}
Recall that $\mathbf X\in \partial g^*_\mu(\mathbf Y)$ if and only if $Y=\nabla g_\mu(\mathbf X)$. The equation $\nabla g_\mu(\mathbf X)= \mathbf Y$ can be written as
\begin{align*}
&\frac{1+\lambda}{\mu}\left[m\mathbf X - \mathbf E\mathbf A\right] + \frac{2}{\mu}\widetilde{\mathbf E}\mathbf X =\mathbf Y\\
&(1+\lambda)\left[m\mathbf X - \mathbf E\mathbf A\right] + 2\widetilde{\mathbf E}\mathbf X =\mu \mathbf Y\\
&\left(m(1+\lambda)\mathbb{I}+ 2\widetilde{\mathbf E}\right)\mathbf X =  (1+\lambda)\mathbf E \mathbf A +\mu \mathbf Y.
\end{align*}
Solving this equation using Lemma \ref{lm51}(ii) yields
\begin{equation}\label{fg}
\mathbf X = \left(\alpha\mathbb{I}+\beta \mathbf E\right)\Big((1+\lambda)\mathbf E \mathbf A + \mu \mathbf Y\Big),
\end{equation}
where $\alpha$ and $\beta$ are given in \eqref{ab}. It follows that $\partial g^*_\mu(\mathbf Y)$ is a singleton for every $\mathbf Y\in \R^{k\times n}$, and so $g^*_\mu$ is continuously differentiable and $\nabla g^*_\mu(\mathbf Y)$ is given by the expression on the right-hand side of \eqref{fg}; see \cite[Theorem 3.3]{BMN14B}. $\h$

To implement the DCA, it remains to find a subgradient of $h_\mu$. From their representations, one can see that $h_{1\mu}$ and  $h_{2\mu}$ are differentiable. Their respective  subgradients coincides with their gradients, that can be computed by  the \emph{partial derivatives} with respect to $x^1, \ldots, x^{k+1}$ given by
\begin{equation}\label{pd1}
\frac{\partial h_{1\mu}}{\partial x^\ell }(\mathbf X) \; =\; (1+\lambda) \sum_{i=1}^{m}\left[\frac{x^\ell - a^i}{\mu}\; - P\left(\frac{x^\ell - a^i}{\mu}; F \right)\right]\;\; \text{for} \;\ell = 1, \ldots, k+1.
\end{equation}
Thus, $\nabla h_{1\mu}(\mathbf X))$ is the $(k+1)\times n$ matrix $\mathbf H_{1}$ whose $\ell^{th}$ row is $\frac{\partial h_{21\mu}}{\partial x^\ell }(\mathbf X)$.

Similarly,
\begin{equation}\label{pd2}
\frac{\partial h_{2\mu}}{\partial x^\ell }(\mathbf X) \; =\; 2 \sum_{j=1}^{k+1}\left[\frac{x^\ell - x^j}{\mu}\; - P\left(\frac{x^\ell - x^j}{\mu}; F \right)\right]\;\; \text{for} \;\ell = 1, \ldots, k+1.
\end{equation}
Hence, $\nabla h_{2\mu}(\mathbf X)$ is the $(k+1)\times n$ matrix $\mathbf H_{4}$ whose $\ell^{th}$ row is $\frac{\partial h_{2\mu}}{\partial x^\ell }(\mathbf X)$.

The procedures for computing a subgradient of  $h_i$ for $i=3, 4, 5$ are similar to those from the previous section. Therefore, we are ready to give a new DCA-based algorithm for the  bilevel hierarchical clustering problem in Model II.

\begin{algorithm}
  \caption{$\;$ {\bf Model~II}}
  \begin{algorithmic}[1]
    \State \textbf{Input:} $ \mathbf A, \mathbf X_0, \lambda_0, \mu_0,\sigma_1, \sigma_2, \epsilon, N \in \mathbb{N}$.
    \State {\bf while} \mbox{stopping criteria} ($\lambda$, $\mu$, $\epsilon$) = false {\bf do}
    \State \qquad $\alpha:=\frac{1}{m(\lambda+1)+2(k+1)}$
	\State \qquad $\beta:=\frac{2}{m(\lambda+1)[m(\lambda+1)+2(k+1)]}$
    \State \qquad {\bf for} $k = 1, \ldots, N\;\;$ {\bf do}
    \State \qquad \qquad Find $\mathbf Y_k\in \partial h_{\mu}(\mathbf X_{k-1})$
    \State \qquad \qquad $\mathbf X_k = \left(\alpha\mathbb{I}+\beta \mathbf E\right)\Big((1+\lambda)\mathbf E \mathbf A + \mu \mathbf Y_k\Big)$
    \State \qquad {\bf end for}
    \State \qquad update $\lambda \;\;\text{and} \; \;\mu$
    \State {\bf end while}
    \State \textbf{Output:} $X_N$.
  \end{algorithmic}
\end{algorithm}

\begin{Example}{\rm ($\ell^2-$clustering with Algorithm 3). In this example, we consider the  hierarchical clustering problem in Model II for the case where $F$ is the Euclidean closed unit ball in $\R^n$. To implement Algorithm 3, it remains to find a subgradient $\mathbf Y\in \partial h_\mu(\mathbf X)$. Recall that
\begin{equation*}
h_\mu(\mathbf X)=h_{1\mu}(\mathbf X)+h_{2\mu}(\mathbf X)+h_3(\mathbf X)+h_4(\mathbf X)+h_5(\mathbf X) \; \mbox{\rm for }\mathbf X\in \R^{(k+1)\times n}.
\end{equation*}
The functions $h_{1\mu}$ and $h_{2\mu}$ are continuously differentiable. The gradients $\nabla h_{1\mu}(\mathbf X)$ and $\nabla h_{2\mu}(\mathbf X)$ can be determined by their partial derivatives from \eqref{pd1} and \eqref{pd2}, respectively. We can find subgradients $\mathbf Y_3\in \partial h_3(\mathbf X)$ and $Y_4\in \partial h_4(\mathbf X)$ by the procedure developed in Example \ref{ex1}. Now, we focus on finding a subgradient $\mathbf Y_5\in \partial h_5(\mathbf X)$. In this case,
\begin{equation*}
h_{5}(\mathbf X):= \max_{t=1, \ldots, k+1}\sum_{\substack{\ell=1\\ \ell\neq t}}^{k+1}\sum_{j=1}^{k+1}\|x^\ell - x^j\|=\max_{t=1, \ldots, k+1}\left(\sum_{\ell=1}^{k+1}\sum_{j=1}^{k+1}\|x^\ell - x^j\|-\sum_{j=1}^{k+1}\|x^t-x^j\|\right).
\end{equation*}
To find such a subgradient, we will apply the subdifferential sum rule  and maximum rule. Choose an index $t^*$ such that
\begin{equation*}
\max_{t=1, \ldots, k+1}\left(\sum_{\ell=1}^{k+1}\sum_{j=1}^{k+1}\|x^\ell - x^j\|-\sum_{j=1}^{k+1}\|x^t-x^j\|\right)=\sum_{\ell=1}^{k+1}\sum_{j=1}^{k+1}\|x^\ell - x^j\|-\sum_{j=1}^{k+1}\|x^{t^*}-x^j\|.
\end{equation*}
Define
\begin{equation*}
v_{\ell j}:=\begin{cases}
      \frac{x^\ell-x^j}{\|x^\ell-x^j\|} &\mbox{\rm if}\;  x^\ell\neq x^j, \\
            0 & \mbox{\rm otherwise}.
   \end{cases}
   \end{equation*}
Then $\mathbf Y_5$ can be determined by the $(k+1)\times n$ matrix whose $\ell^{\rm th}$ row   is given by
\begin{equation*}
Y_\ell:=2\sum_{j=1}^{k+1}v_{\ell j}-v_{\ell t^*}\; \mbox{\rm for }\ell=1, \ldots, k+1.
\end{equation*}
}\end{Example}

By the procedure developed in Example \ref{ex2} with the use of a \emph{signed matrix}, we can similarly provide another example for hierarchical clustering for Model II in the case where $F$ is the closed unit box in $\R^n$. The detail is left for the reader.

\section{Numerical Experiments}
\setcounter{equation}{0}

We conducted our numerical experiments on a MacBook Pro with 2.2 GHz Intel Core i7 Processor, 16 GB 1600 MHz DDR3 Memory. Even though the two continuous optimization formulations we consider are nonsmooth and nonconvex, the Nesterov smoothing technique allowed us to design two implementable DCA-based algorithms.

\begin{figure}[H]
\begin{subfigure}{.50\textwidth}
  \centering
  \includegraphics[width=0.75\textwidth]{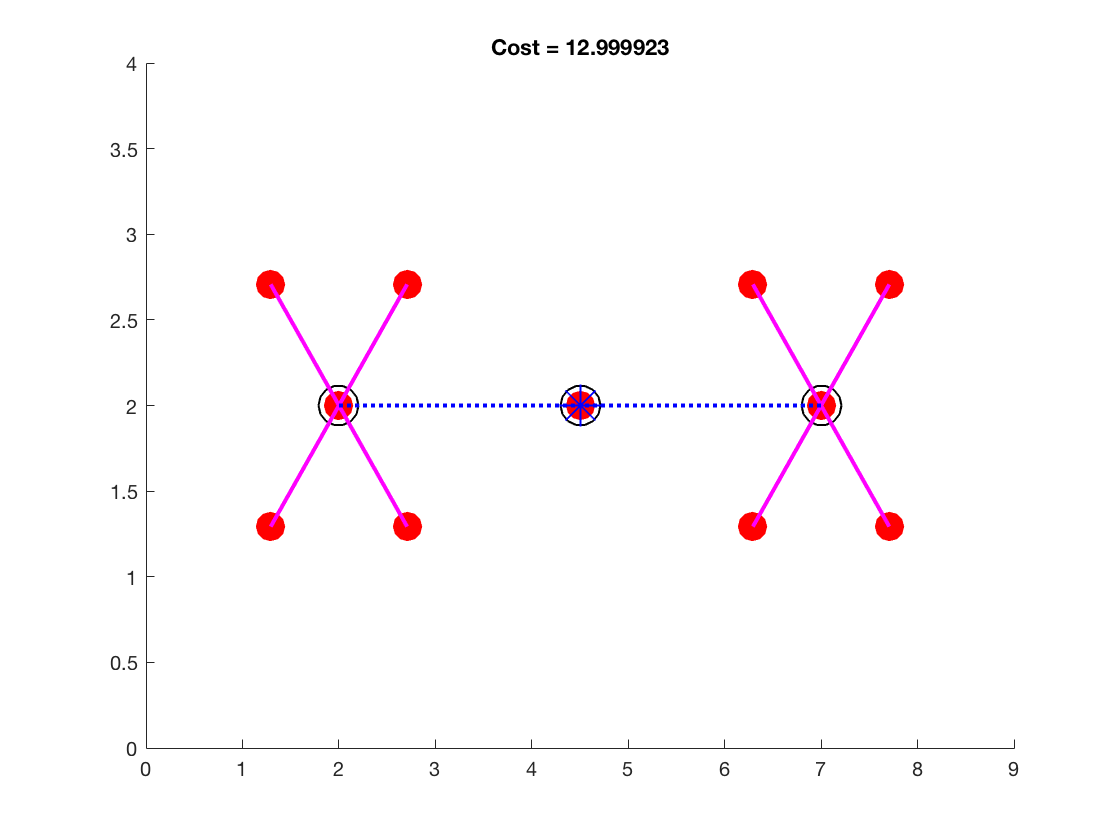}
  \caption{Artificial Test Dataset for  Model~I}
  \label{fig:DS18}
\end{subfigure}%
\begin{subfigure}{.50\textwidth}
  \centering
\includegraphics[width=0.75\textwidth]{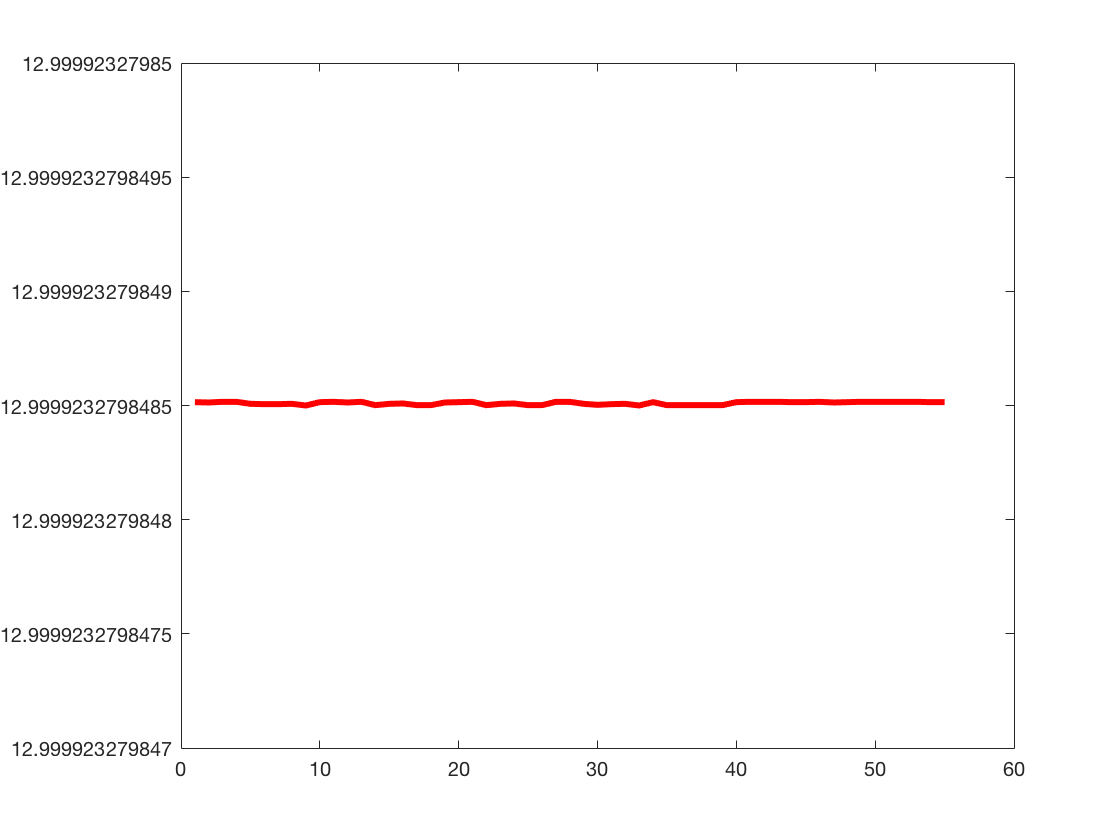}
  \caption{100\% convergence to a global optimal solution}
  \label{fig:DS1002}
\end{subfigure}
\caption{Performance of Algorithm 2.}
\label{fig:alg1}
\end{figure}

For the implementation of the algorithms, we wrote the codes in MATLAB. Since our algorithms are adaptations of the DCA,  there is no guarantee that our algorithms  converge to a global optimal solution. However, for the artificial test dataset we created to test the performance of Algorithm 2 with 11 nodes, 2 clearly identifiable cluster centers, and a total center (see Figure \ref{fig:alg1}), the  algorithm converges 100\% of the time to a global optimal solution for all 55 different pairs of starting centers selected from the 11 points, i.e., $ \binom{11}{2} = 55 $.

\begin{figure}[H]
\begin{subfigure}{.50\textwidth}
  \centering
  \includegraphics[width=0.75\textwidth]{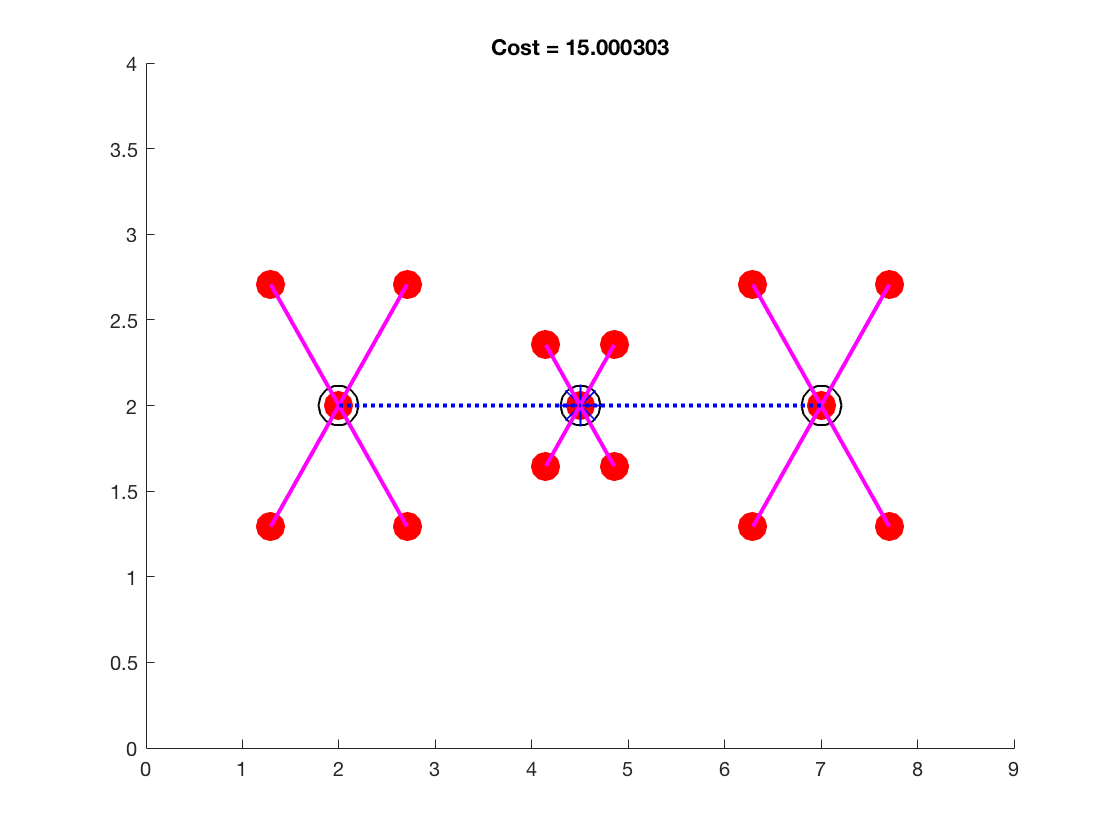}
  \caption{Artificial test dataset for Model~II}
  \label{fig:DS18}
\end{subfigure}%
\begin{subfigure}{.50\textwidth}
  \centering
\includegraphics[width=0.75\textwidth]{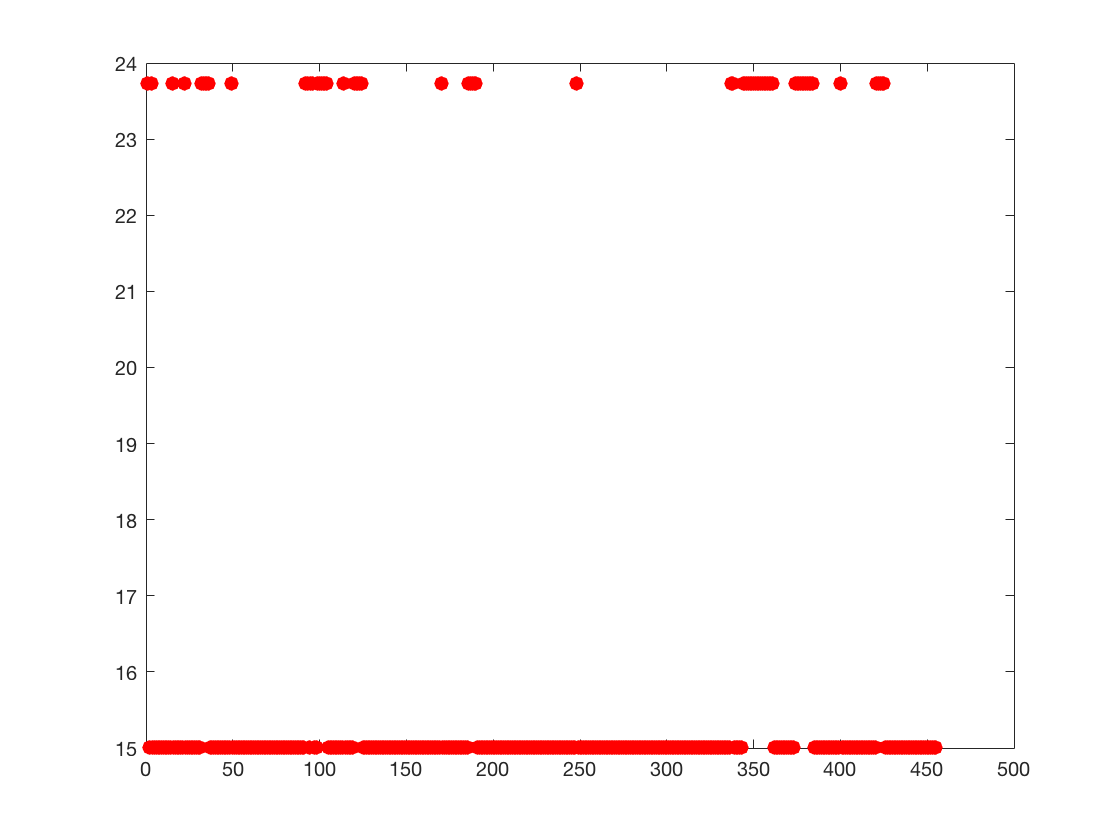}
  \caption{85\% convergence to a global optimal solution}
  \label{fig:DS1002}
\end{subfigure}
\caption{Performance of Algorithm 3 on the Test Data Set.}
\label{fig:alg2}
\end{figure}

On the other hand, for the artificial test dataset we created to test the performance of Algorithm 2 with 15 nodes, 2 clearly identifiable cluster centers, and a total center (see Figure \ref{fig:alg2}), the  algorithm converges to a global optimal solution 85\% of the time, which means that for all 455 different starting centers selected from the 15 points, i.e., $ \binom{15}{3} = 455 $, the algorithm converges to a global optimal solution 85\% of the time.

\begin{figure}[H]
\begin{subfigure}{.50\textwidth}
  \centering
  \includegraphics[width=0.95\textwidth]{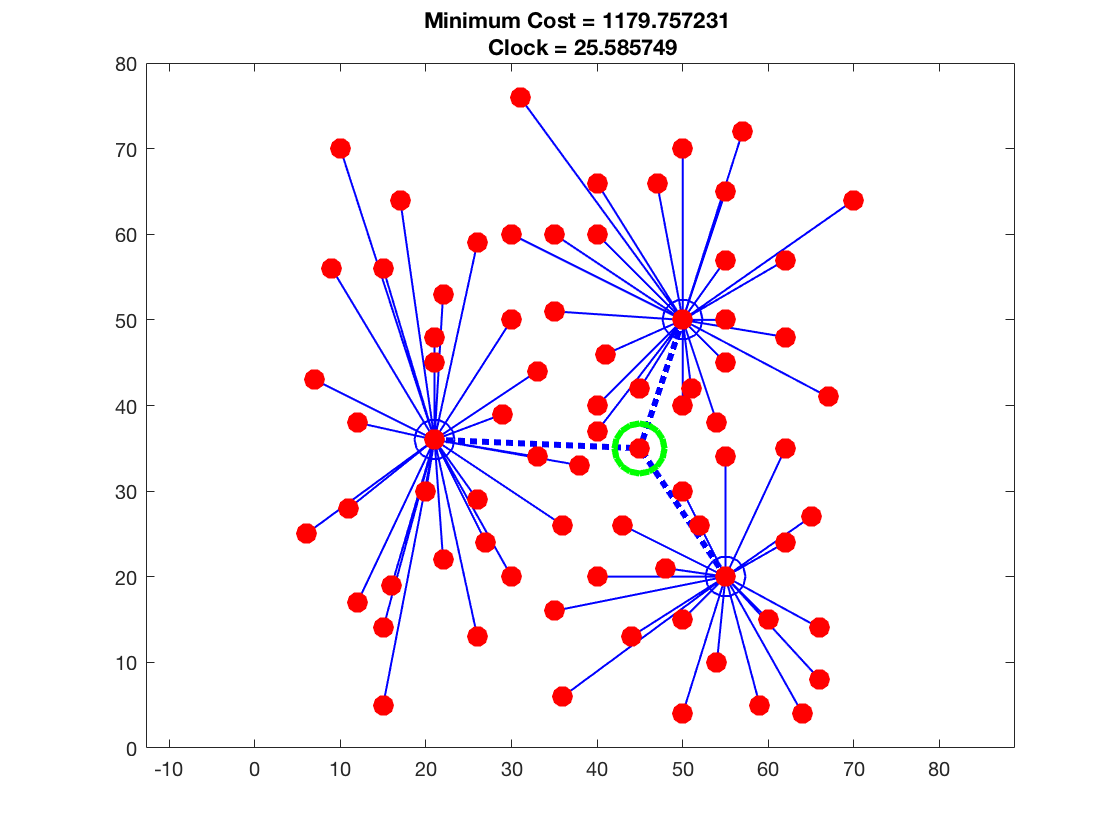}
  \caption{Model I on EIL76}
  \label{fig:DS18}
\end{subfigure}%
\begin{subfigure}{.50\textwidth}
  \centering
\includegraphics[width=0.95\textwidth]{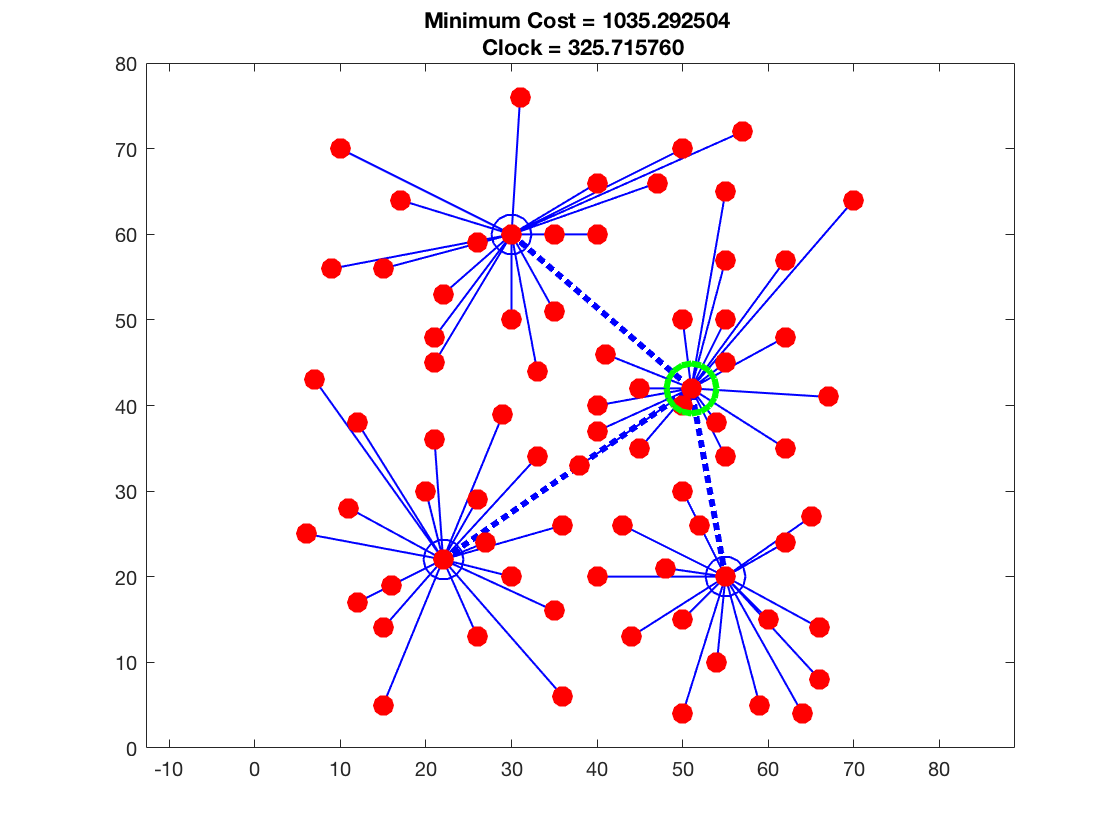}
  \caption{Model II on EIL76}
  \label{fig:DS1002}
\end{subfigure}
\caption{Optimal Solutions for Model I and Model II on EIL76.}
\label{fig:opt}
\end{figure}

Further numerical experiments were performed on the dataset EIL76 (The 76 City Problem) taken from the Traveling Salesman Problem Library \cite{EIL-PR}. For instance, Figures \ref{fig:opt}(a) and \ref{fig:opt}(b) show optimal solutions for Model~I and Model~II, respectively, for three cluster centers and a total center. The optimal solutions were calculated by the brute-force search method in which we exhaustively generated all the four possible candidates, 3 cluster centers and 1 total center, and then computed the corresponding cost to take the minimum. In this case, we have $\binom{76}{3} = 70,300$ combinations for Model~I and $\binom{76}{4} = 1,282,975$ combinations of cluster centers and a total center to check for Model~II. For instance, the optimal value for Model I tested on EIL76 with 3 cluster centers and 1 total center is 1179.76, while for Model II with 3 cluster centers and 1 total center, it is 1035.29.

In the two MATLAB codes we wrote to implement the two algorithms, we updated the penalty parameter $\lambda$ and the smoothing parameter $\mu$ in every iteration by the relations $\lambda_{i+1}=\sigma_1\lambda_{i},\;\sigma_1>1$, and $\mu_{i+1}=\sigma_2\mu_{i},\;\sigma_2\in (0,1)$, respectively. The two parameters were updated until $\mu < 10^{-6}$.

For the choice of the starting centers, we used three different methods:

\begin{itemize}

\item {\bf Random.} We used the ``datasample'' (a MATLAB built in function) to randomly select starting centers from the existing nodes without replacement.

\item {\bf K-means clustering.} We used the ``kmeans'' (a MATLAB built in function) to partition the nodes into $k$ clusters first, and then we selected the $k$ cluster centroid locations as starting centers.

\item {\bf C++ implementation} We implemented the model 1 and model 2 algorithms in C++ and used uniform random numbers generator to generate starting centers. The code was developed using Armadillo library and run on a computer having 20 Intel(R) Xeon(R) CPU E5-2640 v4 @ 2.40GHz cores and 250 GB RAM.

%\item {\bf Genetic Algorithm.}

\end{itemize}

\begin{table}[H]
   \centering
   $\mu_0 = 16,\;\lambda_0 = 0.01,\; \sigma_1  = 160,\;\sigma_2 = 0.5$ %
  \begin{tabular}{ c c c | c  c | c c | c  c c  }
\hline
		&COST1		&COST2			&Time1	&Time2	&Iter1	&Iter2	&k		&m		&n	 \\ \hline
EIL76	&1194.29		&1048.41			&8.04	&10.55	&1058	&1361	&3		&76		&2	 \\
EIL76	&1201.97		&1048.62			&6.84	&7.84	&918		&1006	&3		&76		&2	 \\
EIL76	&1179.76		&1041.53			&7.31	&10.93	&986		&1413	&3		&76		&2	 \\
EIL76	&1181.02		&1057.87			&7.99	&7.50	&1030	&929		&3		&76		&2	 \\
EIL76	&1208.39		&1057.87			&6.40	&7.57	&832		&925		&3		&76		&2	 \\
EIL76	&1179.76		&1057.87			&8.16	&6.77	&1030	&876		&3		&76		&2	 \\
EIL76	&1194.29		&1091.57			&7.89	&6.81	&1056	&881		&3		&76		&2	 \\
EIL76	&1179.76		&1057.87			&7.36	&7.19	&987		&927		&3		&76		&2	 \\
EIL76	&1204.35		&1119.50			&9.97	&9.62	&1337	&1238	&3		&76		&2	 \\
EIL76	&1201.97		&1054.90			&6.98	&6.42	&928		&820		&3		&76		&2	 \\ \hline
\end{tabular}
   \caption{Starting centers  selected randomly, MATLAB code.}
   \label{tab:DS18R}
\end{table}

\begin{table}[H]
   \centering
   $\mu_0 = 16,\;\lambda_0 = 0.01,\; \sigma_1  = 160,\;\sigma_2 = 0.5$ %
  \begin{tabular}{ c c c | c  c | c c | c  c c  }
\hline
		&COST1		&COST2		&Time1	&Time2	&Iter1	&Iter2	&k	&m	&n	 \\ \hline
EIL76	&1204.35	&1059.01	&9.91	&6.62	&1320	&853	&3	&76	&2	 \\
EIL76	&1179.76	&1045.90	&7.23	&9.29	&969	&1195	&3	&76	&2	 \\
EIL76	&1194.29	&1049.53	&7.84	&5.75	&1051	&738	&3	&76	&2	 \\
EIL76	&1179.76	&1059.01	&7.47	&6.61	&994	&853	&3	&76	&2	 \\
EIL76	&1204.35	&1059.01	&9.89	&6.59	&1320	&846	&3	&76	&2	 \\
EIL76	&1179.76	&1059.01	&7.42	&6.64	&994	&853	&3	&76	&2	 \\
EIL76	&1181.02	&1041.29	&7.21	&6.18	&965	&797	&3	&76	&2	 \\
EIL76	&1201.97	&1059.01	&6.99	&6.57	&931	&846	&3	&76	&2	 \\
EIL76	&1181.02	&1059.01	&7.39	&6.62	&988	&853	&3	&76	&2	 \\
EIL76	&1201.97	&1048.62	&6.49	&6.67	&870	&860	&3	&76	&2	 \\ \hline
\end{tabular}
   \caption{Starting centers  selected by the k-means, MATLAB code.}
   \label{tab:DS18R}
\end{table}

\begin{table}[H]
   \centering
   $\mu_0 = 16,\;\lambda_0 = 0.01,\; \sigma_1  = 160,\;\sigma_2 = 0.5$ %
  \begin{tabular}{ c c c | c  c | c c | c  c c  }
\hline
		&COST1		&COST2			&Iter1	&Iter2&Time1	&Time2	&k		&m		&n	 \\ \hline
EIL76 & 1224.04 & 1064.91 & 952 & 829 & 0.09 & 0.05 & 3 & 76 & 2 \\
EIL76 & 1195.55 & 1053.38 & 1051 & 874 & 0.07 & 0.05 & 3 & 76 & 2 \\
EIL76 & 1206.92 & 1041.52 & 1045 & 1091 & 0.07 & 0.07 & 3 & 76 & 2 \\
EIL76 & 1206.92 & 1057.86 & 1008 & 855 & 0.06 & 0.06 & 3 & 76 & 2 \\
EIL76 & 1215.56 & 1065.79 & 1165 & 887 & 0.07 & 0.05 & 3 & 76 & 2 \\
EIL76 & 1218.48 & 1057.86 & 1263 & 829 & 0.07 & 0.04 & 3 & 76 & 2 \\
EIL76 & 1197.42 & 1067.6 & 988 & 884 & 0.04 & 0.04 & 3 & 76 & 2 \\
EIL76 & 1206.92 & 1048.6 & 1045 & 1020 & 0.05 & 0.04 & 3 & 76 & 2 \\
EIL76 & 1215.56 & 1057.86 & 1148 & 843 & 0.05 & 0.04 & 3 & 76 & 2 \\
EIL76 & 1215.56 & 1165.62 & 1206 & 920 & 0.05 & 0.04 & 3 & 76 & 2 \\

\hline
\end{tabular}
   \caption{Starting centers  selected randomly, C++ code}
   \label{tab:DS18R}
\end{table}

\begin{table}[H]
   \centering
   $\mu_0 = 16,\;\lambda_0 = 0.01,\; \sigma_1  = 160,\;\sigma_2 = 0.5$ %
  \begin{tabular}{ c c c | c  c | c c | c  c c  }
\hline
		&COST1		&COST2			&Iter1	&Iter2&Time1	&Time2	&k		&m		&n	 \\ \hline
1002C & 2.56341e+06 & 2.24537e+06 & 1023 & 1023 & 1.31 & 1 & 6 & 1002 & 2 \\
1002C & 2.16241e+06 & 1.79317e+06 & 1023 & 1023 & 1.09 & 1 & 6 & 1002 & 2 \\
1002C & 2.55508e+06 & 2.25252e+06 & 1023 & 1023 & 1.1 & 0.99 & 6 & 1002 & 2 \\
1002C & 2.29283e+06 & 2.12459e+06 & 1023 & 1023 & 1.1 & 0.99 & 6 & 1002 & 2 \\
1002C & 2.28579e+06 & 2.02933e+06 & 1023 & 1023 & 1.1 & 1 & 6 & 1002 & 2 \\
1002C & 2.02867e+06 & 1.84531e+06 & 1023 & 1023 & 1.1 & 0.99 & 6 & 1002 & 2 \\
1002C & 2.49236e+06 & 2.43734e+06 & 1023 & 1023 & 1.1 & 0.99 & 6 & 1002 & 2 \\
1002C & 3.02324e+06 & 2.42825e+06 & 1023 & 1023 & 1.1 & 0.99 & 6 & 1002 & 2 \\
1002C & 2.33796e+06 & 2.1374e+06 & 1023 & 1023 & 1.1 & 1 & 6 & 1002 & 2 \\
1002C & 2.37677e+06 & 1.85446e+06 & 1023 & 1023 & 1.09 & 1 & 6 & 1002 & 2 \\

\hline
\end{tabular}
   \caption{Starting centers  selected randomly, C++ code}
   \label{tab:DS18R}
\end{table}

\begin{table}[H]
   \centering
   $\mu_0 = 16,\;\lambda_0 = 0.01,\; \sigma_1  = 160,\;\sigma_2 = 0.5$ %
  \begin{tabular}{ c c c | c  c | c c | c  c c  }
\hline
		&COST1		&COST2				&Iter1	&Iter2&Time1	&Time2	&k		&m		&n	 \\ \hline
10000RND & 1.94933e+07 & 1.8097e+07 & 1023 & 1023 & 11.36 & 10.1 & 6 & 10000 & 2 \\
10000RND & 2.44543e+07 & 2.07372e+07 & 1023 & 1023 & 11.15 & 10.1 & 6 & 10000 & 2 \\
10000RND & 2.36188e+07 & 1.90255e+07 & 1023 & 1023 & 11.18 & 10.07 & 6 & 10000 & 2 \\
10000RND & 2.13395e+07 & 1.81326e+07 & 1023 & 1023 & 11.16 & 10.09 & 6 & 10000 & 2 \\
10000RND & 1.97625e+07 & 1.74163e+07 & 1023 & 1023 & 11.17 & 10.09 & 6 & 10000 & 2 \\
10000RND & 1.9848e+07 & 1.79588e+07 & 1023 & 1023 & 11.18 & 10.11 & 6 & 10000 & 2 \\
10000RND & 2.4502e+07 & 2.0164e+07 & 1023 & 1023 & 11.17 & 10.08 & 6 & 10000 & 2 \\
10000RND & 2.38836e+07 & 2.09025e+07 & 1023 & 1023 & 11.16 & 10.09 & 6 & 10000 & 2 \\
10000RND & 1.81975e+07 & 1.68355e+07 & 1023 & 1023 & 11.17 & 10.09 & 6 & 10000 & 2 \\
10000RND & 2.05324e+07 & 1.68926e+07 & 1023 & 1023 & 11.16 & 10.1 & 6 & 10000 & 2 \\

\hline
\end{tabular}
   \caption{Starting centers  selected randomly, C++ code, 10000 u. randomly distributed points}
   \label{tab:DS18R}
\end{table}

\begin{table}[H]
   \centering
   $\mu_0 = 16,\;\lambda_0 = 0.01,\; \sigma_1  = 160,\;\sigma_2 = 0.5$ %
  \begin{tabular}{ c c c | c  c | c c | c  c c  }
\hline
		&COST1		&COST2				&Iter1	&Iter2	&Time1	&Time2 &k		&m		&n	 \\ \hline
10000RND & 5.17176e+06 & 5.10097e+06 & 1023 & 1023 & 218.72 & 166.85 & 100 & 10000 & 2 \\
10000RND & 5.32321e+06 & 5.20111e+06 & 1023 & 1023 & 218.1 & 164.76 & 100 & 10000 & 2 \\
10000RND & 5.32893e+06 & 5.21018e+06 & 1023 & 1023 & 215.79 & 166.91 & 100 & 10000 & 2 \\
10000RND & 5.45463e+06 & 5.34531e+06 & 1023 & 1023 & 217.58 & 166.92 & 100 & 10000 & 2 \\
10000RND & 5.59697e+06 & 5.42149e+06 & 1023 & 1023 & 217.25 & 164.93 & 100 & 10000 & 2 \\
10000RND & 5.57053e+06 & 5.39613e+06 & 1023 & 1023 & 215.23 & 169.07 & 100 & 10000 & 2 \\
10000RND & 5.67843e+06 & 5.55442e+06 & 1023 & 1023 & 217.15 & 166.78 & 100 & 10000 & 2 \\
10000RND & 5.7148e+06 & 5.57767e+06 & 1023 & 1023 & 215.6 & 165.05 & 100 & 10000 & 2 \\
10000RND & 5.37335e+06 & 5.28977e+06 & 1023 & 1023 & 219.63 & 164.81 & 100 & 10000 & 2 \\
10000RND & 5.73865e+06 & 5.61554e+06 & 1023 & 1023 & 217.23 & 166.87 & 100 & 10000 & 2 \\

\hline
\end{tabular}
   \caption{Starting centers  selected randomly, C++ code, 10000 u. randomly distributed points}
   \label{tab:DS18R}
\end{table}

\begin{table}[H]
   \centering
   $\mu_0 = 16,\;\lambda_0 = 0.01,\; \sigma_1  = 160,\;\sigma_2 = 0.5$ %
  \begin{tabular}{ c c c | c  c | c c | c  c c  }
\hline
		&COST1		&COST2				&Iter1	&Iter2	&Time1	&Time2 &k		&m		&n	 \\ \hline
10000RND3D & 2.83948e+07 & 2.63213e+07 & 1023 & 1023 & 24.79 & 19.71 & 10 & 10000 & 3 \\
10000RND3D & 2.74404e+07 & 2.65681e+07 & 1023 & 1023 & 24.6 & 19.7 & 10 & 10000 & 3 \\
10000RND3D & 2.9869e+07 & 2.85641e+07 & 1023 & 1023 & 24.59 & 19.7 & 10 & 10000 & 3 \\
10000RND3D & 3.44097e+07 & 3.07609e+07 & 1023 & 1023 & 24.6 & 19.7 & 10 & 10000 & 3 \\
10000RND3D & 3.05076e+07 & 2.89047e+07 & 1023 & 1023 & 24.6 & 19.7 & 10 & 10000 & 3 \\
10000RND3D & 2.72841e+07 & 2.61452e+07 & 1023 & 1023 & 24.61 & 19.7 & 10 & 10000 & 3 \\
10000RND3D & 2.94171e+07 & 2.81767e+07 & 1023 & 1023 & 24.6 & 22.25 & 10 & 10000 & 3 \\
10000RND3D & 3.15467e+07 & 2.72963e+07 & 1023 & 1023 & 24.6 & 19.69 & 10 & 10000 & 3 \\
10000RND3D & 2.78719e+07 & 2.64644e+07 & 1023 & 1023 & 24.61 & 21.48 & 10 & 10000 & 3 \\
10000RND3D & 2.80267e+07 & 2.64164e+07 & 1023 & 1023 & 24.59 & 19.7 & 10 & 10000 & 3 \\

\hline
\end{tabular}
   \caption{Starting centers  selected randomly, C++ code, 10000 u. randomly distributed points, 3 dimensions}
   \label{tab:DS18R}
\end{table}

\begin{table}[H]
   \centering
   $\mu_0 = 16,\;\lambda_0 = 0.01,\; \sigma_1  = 160,\;\sigma_2 = 0.5$ %
  \begin{tabular}{ c c c | c  c | c c | c  c c  }
\hline
		&COST1		&COST2				&Iter1	&Iter2	&Time1	&Time2 &k		&m		&n	 \\ \hline
1000RND6D & 5.68343e+06 & 5.49334e+06 & 1023 & 1023 & 2.72 & 2.16 & 10 & 1000 & 6 \\
1000RND6D & 6.15169e+06 & 5.94648e+06 & 1023 & 1023 & 2.5 & 2.15 & 10 & 1000 & 6 \\
1000RND6D & 5.95467e+06 & 5.87668e+06 & 1023 & 1023 & 2.51 & 2.15 & 10 & 1000 & 6 \\
1000RND6D & 5.848e+06 & 5.67641e+06 & 1023 & 1023 & 2.5 & 2.16 & 10 & 1000 & 6 \\
1000RND6D & 5.82286e+06 & 5.73382e+06 & 1023 & 1023 & 2.5 & 2.15 & 10 & 1000 & 6 \\
1000RND6D & 5.81637e+06 & 5.49823e+06 & 1023 & 1023 & 2.51 & 2.15 & 10 & 1000 & 6 \\
1000RND6D & 6.00205e+06 & 5.84304e+06 & 1023 & 1023 & 2.5 & 2.15 & 10 & 1000 & 6 \\
1000RND6D & 5.9963e+06 & 5.86284e+06 & 1023 & 1023 & 2.5 & 2.17 & 10 & 1000 & 6 \\
1000RND6D & 6.16517e+06 & 6.03364e+06 & 1023 & 1023 & 2.5 & 2.14 & 10 & 1000 & 6 \\
1000RND6D & 5.71309e+06 & 5.60686e+06 & 1023 & 1023 & 2.51 & 2.15 & 10 & 1000 & 6 \\

\hline
\end{tabular}
   \caption{Starting centers  selected randomly, C++ code, 1000 u. randomly distributed points in 6 dimensions}
   \label{tab:DS18R}
\end{table}

\begin{table}[H]
   \centering
   $\mu_0 = 16,\;\lambda_0 = 0.01,\; \sigma_1  = 160,\;\sigma_2 = 0.5$ %
  \begin{tabular}{ c c c | c  c | c c | c  c c  }
\hline
		&COST1		&COST2				&Iter1	&Iter2	&Time1	&Time2 &k		&m		&n	 \\ \hline
100000RND2D & 1.40282e+08 & 1.33498e+08 & 1023 & 1023 & 198.3 & 165.35 & 10 & 100000 & 2 \\
100000RND2D & 1.83297e+08 & 1.54512e+08 & 1023 & 1023 & 197.06 & 168.74 & 10 & 100000 & 2 \\
100000RND2D & 1.5134e+08 & 1.41451e+08 & 1023 & 1023 & 198.74 & 165.34 & 10 & 100000 & 2 \\
100000RND2D & 1.59333e+08 & 1.4203e+08 & 1023 & 1023 & 199.65 & 164.96 & 10 & 100000 & 2 \\
100000RND2D & 1.53366e+08 & 1.35764e+08 & 1023 & 1023 & 199.08 & 167.07 & 10 & 100000 & 2 \\
100000RND2D & 1.55465e+08 & 1.45342e+08 & 1023 & 1023 & 199.99 & 166.82 & 10 & 100000 & 2 \\
100000RND2D & 1.39211e+08 & 1.32843e+08 & 1023 & 1023 & 197.37 & 165.71 & 10 & 100000 & 2 \\
100000RND2D & 1.60153e+08 & 1.4911e+08 & 1023 & 1023 & 199.78 & 167.28 & 10 & 100000 & 2 \\
100000RND2D & 1.52469e+08 & 1.38242e+08 & 1023 & 1023 & 200.14 & 167.13 & 10 & 100000 & 2 \\
100000RND2D & 1.46638e+08 & 1.38241e+08 & 1023 & 1023 & 197.63 & 165.07 & 10 & 100000 & 2 \\

\hline
\end{tabular}
   \caption{Starting centers  selected randomly, C++ code, 100000 u. randomly distributed points in 2 dimensions}
   \label{tab:DS18R}
\end{table}

\section{Conclusion and Future Research}

In this study, we presented two DCA-based algorithms for solving two different bilevel hierarchical clustering problems where the similarity(dissimilarity) measure between two data points (nodes) is given by generalized distances. As special cases of generalized distances, we  provided two detailed examples for the $\ell^1$ and $\ell^2$ norms. We implemented the algorithms with MATLAB and C++ and tested them on different datasets of various sizes and dimensions. We expect that our method used in this paper for solving bilevel hierarchical clustering problems are applicable to solving other nonsmooth nonconvex optimization problems. 
%%%%%%%%%%%%%%%%%%%%%%%%%%%%%%%%%


\small
\begin{thebibliography}{99}

\bibitem{abt07} An, L.T.H., Belghiti, M.T., Tao, P.D.:  A new efficient algorithm based on DC programming and DCA for clustering. J. Glob. Optim., \textbf{27},  503--608 (2007).

\bibitem{amt14} An, L.T.H., Minh, L.H., Tao, P.D.: New and efficient DCA based algorithms for minimum sum-of-squares clustering, Pattern Recognition, \textbf{47}, 388--401(2014).

\bibitem{am07}  An, L.T.H., Minh, L.H.: Optimization based DC programming and DCA for hierarchical clustering. European J. Oper. Res. \textbf{183}, 1067--1085 (2007).

\bibitem{at97} An, L.T.H., Tao, P.D.:  Convex analysis approach to D.C. programming: Theory, algorithms and applications. Acta Math. Vietnam. \textbf{22}, 289--355 (1997).

\bibitem{B99} Bagirov, A.: Derivative-free methods for unconstrained nonsmooth optimization and its numerical analysis. Investigacao Operacional. \textbf{19}, 75--93 (1999).

\bibitem{bjor03} Bagirov, A., Jia, L., Ouveysi, I., Rubinov, A.M.: Optimization based clustering algorithms in Multicast group hierarchies, in: Proceedings of the Australian Telecommunications, Networks and Applications Conference (ATNAC), Melbourne Australia (published on CD, ISNB 0-646-42229-4) (2003).

\bibitem{btu16} Bagirov, A.,  Taheri, S., Ugon, J.:  Nonsmooth DC programming approach to the minimum sum-of-squares clustering problems. Pattern Recognition. \textbf{53}, 12--24 (2016).

\bibitem{bvx13} Barbosa, G. V., Villas-Boas, S. B., Xavier, A. E.: Solving the Two-level Clustering Problem by Hyperbolic Smoothing Approach, and Design of Multicast Networks, SELECTED PROCEEDINGS, WCTR RIO (2013).

\bibitem{bc11} Bauschke, H.H.,  Combettes, P.L.: Convex Analysis and Monotone Operator Theory in Hilbert Spaces. Springer, New York (2011).

\bibitem{bl06} Borwein, J.M., Lewis, A.S.: Convex Analysis and Nonlinear Optimization, 2nd edition. Springer, New York (2006).

\bibitem{b10} Bo\c t, R.I.:  Conjugate Duality in Convex Optimization.  Springer, Berlin (2010).

\bibitem{h59} Hartman, P.:  On functions representable as a difference of convex functions. Pacific J. Math. \textbf{9}, 707--713 (1959).

\bibitem{hl93} Hiriart-Urruty, J.B.,  Lemar\'echal, C.:  Convex Analysis and Minimization Algorithms I, II. Springer, Berlin (1993).

\bibitem{h85} Hiriart-Urruty, J.B.:  Generalized differentiability, duality and optimization for problems dealing with differences of convex functions. Lecture Note in Economics and Math. Systems. \textbf{256}, 37--70 (1985).

\bibitem{BM06B} Mordukhovich, B.S.:  Variational Analysis and Generalized Differentiation, I: Basic Theory, II: Applications. Springer, Berlin (2006).

\bibitem{BMN14B} Mordukhovich, B.S.,  Nam, N.M.:  An Easy Path to Convex Analysis and Applications. Morgan \& Claypool Publishers, San Rafael, CA (2014).

\bibitem{nars14} Nam, N. M., An, N. T., Rector, R. B., J. Sun, J.: Nonsmooth algorithms and Nesterov's smoothing technique for generalized Fermat-Torricelli problems. SIAM J. Optim. \textbf{24}, 1815--1839 (2014).

\bibitem{nrg17} Nam, N.M., Rector, R.B.,  Giles, D.: Minimizing Differences of Convex Functions with Applications to Facility Location and Clustering. Journal of Optimization Theory and Applications. \textbf{173}, 255--278 (2017).

\bibitem{wondi} Nam, N. M., Geremew, W., Reynolds, S., Tran, T: The Nesterov Smoothing Technique and Minimizing Differences of Convex Functions for Hierarchical Clustering, in press (2017).

\bibitem{ny13} Nesterov, Y.:   Gradient methods for minimizing composite functions. Math. Program. \textbf{140}, 125--161 (2013).

\bibitem{ny05} Nesterov, Y.:  Smooth minimization of non-smooth functions. Math. Program. \textbf{103}, 127--152 (2005).

\bibitem{ny04b} Nesterov, Y.: Introductory lectures on convex optimization. A basic course. Applied Optimization, 87. Kluwer Academic Publishers, Boston, MA (2004).

\bibitem{ob15} Ordin, B.,  Bagirov, A.:  A heuristic algorithm for solving the minimum sum-of-squares clustering problems. Journal of Global Optimization. \textbf{61}, 341--361 (2015).

\bibitem{EIL-PR} Reinelt, G.: TSPLIB: A Traveling Salesman Problem Library.  ORSA Journal of Computing. \textbf{3}, 376--384 (1991).

\bibitem{r70b} Rockafellar, R.T.:  Convex Analysis. Princeton University Press, Princeton, NJ (1970).

\bibitem{r74b} Rockafellar, R.T.:  Conjugate Duality and Optimization. SIAM, Philadelphia, PA (1974).

\bibitem{ta98} Tao, P.D.,  An, L.T.H.: A d.c. optimization algorithm for solving the trust-region subproblem, SIAM J. Optim. \textbf{8}, 476--505 (1998).


\end{thebibliography}
\end{document}